\documentclass[12pt]{amsart}
\usepackage{amsmath,amssymb,amsthm,hyperref,color,euscript}

\newcommand{\ds}{\displaystyle}

\newcommand{\cd}{\cdot}
\newcommand{\cds}{\cdots}

\newcommand{\vsb}{\vspace{2mm}}

\newcommand{\q}{\quad}

\newcommand{\maru}[1]{{\ooalign{\hfil#1\/\hfil\crcr
\raise.167ex\hbox{\mathhexbox20D}}}}

\newcommand{\ruby}[2]{%
 \leavevmode
 \setbox0=\hbox{#1}%
 \setbox1=\hbox{\tiny #2}%
 \ifdim\wd0>\wd1 \dimen0=\wd0 \else \dimen0=\wd1 \fi
 \hbox{%
   \kanjiskip=0pt plus 2fil
   \xkanjiskip=0pt plus 2fil
   \vbox{%
     \hbox to \dimen0{%
       \tiny \hfil#2\hfil}%
     \nointerlineskip
     \hbox to \dimen0{\mathstrut\hfil#1\hfil}}}}

\newcommand{\la}{\langle}
\newcommand{\ra}{\rangle}

\DeclareMathOperator*{\tensor}{\otimes}
\DeclareMathOperator*{\fusion}{\boxtimes}

\newcommand{\Z}{\mathbb{Z}}
\newcommand{\C}{\mathbb{C}}
\newcommand{\R}{\mathbb{R}}

\newcommand{\RM}{\mathrm{RM}}
\newcommand{\calC}{\EuScript{C}}

\newcommand{\Sym}{{\rm Sym}}
\newcommand{\Imm}{{\rm Im}}

\newcommand{\vir}{\mathrm{Vir}}
\newcommand{\aut}{\mathrm{Aut}}
\newcommand{\Aut}{\mathrm{Aut}}
\newcommand{\wt}{\mathrm{wt}}

\newcommand{\om}{\omega}
\newcommand{\be}{\beta}
\newcommand{\al}{\alpha}

\newcommand{\Span}{\mathrm{Span}}

\newcommand{\w}{\omega}

\makeatletter \@addtoreset{equation}{section}

\theoremstyle{plain}
\newtheorem{theorem}{Theorem}[section]
\newtheorem{proposition}[theorem]{Proposition}
\newtheorem{lemma}[theorem]{Lemma}
\newtheorem{corollary}[theorem]{Corollary}

\theoremstyle{definition}
\newtheorem{definition}[theorem]{Definition}

\newtheorem{example}{Example}
\newtheorem{notation}[theorem]{Notation}
\theoremstyle{remark}
\newtheorem{remark}[theorem]{Remark}
\numberwithin{equation}{section}

\title[Frame Stabilizer]{Frame Stabilizers for framed vertex operator algebras associated to lattices having 4-frames}
\author{Ching Hung Lam} %
  \address[C. H. Lam] {Department of Mathematics and National Center for Theoretical Sciences, National Cheng Kung University, Tainan, Taiwan 701}
  \email{chlam@mail.ncku.edu.tw}
\author[H. Shimakura]{Hiroki Shimakura}%
\address[H. Shimakura]{Department of Mathematics,
Aichi University of Education,
1 Hirosawa, Igaya-cho, Kariya-city, Aichi, 448-8542 Japan}%
\email {shima@auecc.aichi-edu.ac.jp}%
\date{}
\thanks{C.\,H. Lam is partially supported by NSC grant
  97-2115-M-006-015-MY3 of Taiwan, R.O.C.}
\thanks{H.\ Shimakura is partially supported by the Grants-in-Aid for Scientific Research No. 20549004.}

\newcommand{\sfr}[2]{\leavevmode\kern-.1em
  \raise.5ex\hbox{\the\scriptfont0 #1}\kern-.1em
  /\kern-.15em\lower.25ex\hbox{\the\scriptfont0 #2}}

\newcommand{\shf}{\sfr{1}{2}}

\newcommand{\stab}{\mathrm{Stab}}
\newcommand{\pstab}{\mathrm{Stab}^{\mathrm{pt}}}

\begin{document}
\begin{abstract}
In this paper, we study certain Virasoro frames for lattice vertex operator algebras and their $\Z_2$-orbifolds using linear codes over $\Z_4$. We also compute the corresponding frame stabilizer from the view point of binary codes and $\Z_4$-codes. As an application, we determine the frame stabilizers of several Virasoro frames of the vertex operator algebra $V_{E_8}$ and the moonshine vertex operator algebra $V^\natural$.
\end{abstract}

\maketitle
\section{Introduction}

A framed vertex operator algebra $V$ is a simple vertex operator
algebra (VOA) which contains a Virasoro frame, a subVOA $T_r$ isomorphic to the
tensor product of $r$-copies of the simple Virasoro VOA $L(\shf,0)$
such that $\mathrm{rank}\ V=\mathrm{rank}\ T_r=r/2$. There are many
important examples such as the moonshine VOA $V^\natural$ and the
Leech lattice VOA. In \cite{DGH}, a basic theory of  framed VOAs was
established. A general structural theory about the automorphism
group and the frame stabilizer, the subgroup which stabilizes $T_r$
setwise, was also included in \cite{DGH,LY}. It was shown in \cite{DGH} that the frame stabilizer of a framed VOA $V=\oplus_{i\in \Z}V_i$ is always a finite group.
Moreover, Miyamoto\,\cite{M3} showed that the full automorphism group $\mathrm{Aut}\, (V)$ of a framed VOA $V$ is also finite if $V_1=0$. Hence, the theory of framed VOA is very
useful in studying certain finite groups such as the Monster.

Lattice VOAs associated to even lattices are basic examples of VOAs \cite{Bo,FLM}.
If an even lattice $L$ of rank $n$ has a $4$-frame, i.e., an orthogonal basis of $\R\otimes_\Z L$ of norm $4$, then the lattice VOA $V_L$ has a natural Virasoro frame $T_{2n}$ (\cite{DMZ}), which is fixed by the involution in $\Aut(V_L)$ lifted from the $-1$-isometry on $L$.
Therefore, if $L$ is unimodular, $T_{2n}$ is also contained in the $\Z_2$-orbifold VOA $\tilde{V}_L=V_L^+\oplus V_L^{T,+}$.
The main purpose of this paper is to determine the frame stabilizer in $\Aut(V_L)$ and $\Aut(\tilde{V}_L)$ of the Virasoro frame $T_{2n}$ associated to a $4$-frame of $L$.

Given a framed VOA $V$, a Virasoro frame $T_r$ determines two binary codes $(C,D)$ of length $r$, called the structure codes. In \cite{LY}, it was shown that the structure codes $C$ and $D$ satisfy certain duality conditions. The pointwise frame stabilizer, the subgroup which fixes $T_r$ pointwise, was also determined.

It is known \cite{GH,LY} that the frame stabilizer naturally acts on $C$ and $D$ while the pointwise frame stabilizer is the subgroup of the frame stabilizer acting on $C$ and $D$ trivially.
The main problem is thus to determine the quotient of the frame stabilizer by the pointwise frame stabilizer, which is isomorphic to a subgroup of $\Aut(C)\cap\Aut(D)$ (\cite{DGH,LY}).

In \cite{GH}, the Virasoro frames for the lattice VOA $V_{E_8}$ were studied. It was shown that there are exactly $5$ Virasoro frames for $V_{E_8}$ up to conjugation. The corresponding frame stabilizers were also computed. It was shown that for any Virasoro frame, the quotient of the frame stabilizer by the pointwise frame stabilizer is equal to the automorphism group of the structure codes. Their method for twisted case $\tilde{V}_{E_8}(\cong V_{E_8})$ associated to the fourth $\Z_4$-code, however, used some special properties of the Lie group $E_8(\C)$, which is the automorphism group of $V_{E_8}$. It may be difficult to generalize to other cases.

\medskip

In this article, we shall study a VOA $V$ isomorphic to the lattice VOA $V_L$ associated to an even lattice $L$ of rank
$n$ having a $4$-frame $F$, and its $\Z_2$-orbifold $\tilde{V}_L$ when $L$ is unimodular.  Since $F^*/F\cong \Z_4^n$, $L/F$ is isomorphic to a
self-orthogonal $\Z_4$-code $\EuScript{C}$. Our approach mainly stresses on the
relationship between the self-orthogonal $\Z_4$-code $\EuScript{C}$
and the structure codes of the corresponding framed VOA. More precisely, we shall determine
the frame stabilizer of $V$ in terms of the automorphism group of the structure code $C$ and the automorphism group of the $\Z_4$-code $\EuScript{C}$. Let $T_{2n}$ be the
Virasoro frame of $V$ associated to $F$. We first study a certain subcode of $C$, which is isomorphic to $d(\Z_2^n)$  for $V= V_L$ and to $d(\mathcal{E}_n)$ for $V= \tilde{V}_L$, where $d:\Z_2^n \to \Z_2^{2n}$ is the ``double'' map defined by $$d(c_1, c_2,\dots,c_n)= (c_1,c_1, c_2, c_2, \dots, c_n,c_n),\qquad  (c_1,\dots,c_n)\in \Z_2^n$$
and $\mathcal{E}_n$ is the subcode of $\Z_2^n$ consisting of all even words (cf. Notation  \ref{Nde} and \ref{NEn}). We shall describe the stabilizer of $d(\Z_2^n)$ for $V= V_L$ (resp. $d(\mathcal{E}_n)$ for $V= \tilde{V}_L$) in the subgroup
$\stab_{\Aut(V)}(T_{2n})/\pstab_{\Aut(V)}(T_{2n})$ of $\Aut(C)$ in terms of the
automorphism group of the $\Z_4$-code $\EuScript{C}$. The main result is the transitivity
of $\stab_{\Aut(V)}(T_{2n})/\pstab_{\Aut(V)}(T_{2n})$ on the set of such subcodes of $C$ (under the assumption that the minimum weight of the binary code $\calC_0$ is greater than or equal to $4$ for $V=\tilde{V}_L$).
In addition, we shall show $\stab_{\Aut(V)}(T_{2n})$ is generated by some automorphisms induced from $\Aut(L)$ and triality automorphisms \cite{FLM}, and determine the index of
$\stab_{\Aut(V)}(T_{2n})/\pstab_{\Aut(V)}(T_{2n})$ in $\Aut(C)$.
Hence, one can determine the frame stabilizer $\stab_{\Aut(V)}(T_{2n})$ in principle by using the automorphism groups $\Aut(C)$ and $\Aut(\EuScript{C})$.

As an application, we shall study some Virasoro frames for the lattice VOA $V_{E_8}$ and certain Virasoro frames of the moonshine VOA $V^\natural$ arisen from a pseudo Golay code over $\Z_4$, that means it reduces to the binary Golay code modulo $2$ and extends to the Leech lattice through the construction $A$. We shall also study the Virasoro frame of $V^\natural$ associated to the standard $4$-frame of the Leech lattice. The corresponding frame stabilizer in each case will also be computed.

\begin{remark}
A binary code $C$ is said to be \textit{indecomposable} if it cannot be written
as a direct sum of two subcodes of shorter length. We believe that any
indecomposable maximal binary code of length $48$ consisting of codewords whose
weights are divisible by $8$ is isomorphic to $\Span_{\Z_2}\{d(B),e(1^{24})\}$
for some doubly even self-dual code $B$ of length $24$ (see Notation \ref{Nde}
for the definitions of the maps $d$ and $e$).
If this is true, then any holomorphic framed VOA of rank $24$ is isomorphic to a
 lattice VOA or its $\Z_2$-orbifold, and any frame is conjugate to the frame
associated to a $4$-frame of the lattice (cf. \cite{La07}). 
Therefore, our method can be applied to the Virasoro frames of most, if not
all, holomorphic framed VOAs of rank $24$.
\end{remark}

\paragraph{\bf Acknowledgements.} The authors thank Professor Masaaki Harada and Professor Akihiro Munemasa for helpful comments on $\Z_4$-codes.
Part of the work was done when the second author was visiting the National Center for Theoretical Sciences, Taiwan on November and December 2008.
He thanks the staff of the center for their help.

\medskip




\paragraph{\bf Notation and terminology}
\begin{center}
\begin{tabular}{ll}
$\langle\cdot,\cdot\rangle$& The inner product in $\Z_2^n$ defined by $\langle X,Y\rangle=|X\cap Y|\mod2$.\\
$(\cdot,\cdot)$& The inner product in $\R^n$.\\
$2^{n}$& An elementary abelian $2$-group of order $2^n$.\\
$2^{n_1+n_2+\dots+n_k}$& A group extension $2^{n_1}.(2^{n_2}.(\cdots.(2^{n_k})\cdots)$.\\
$A_4(\calC)$& The lattice obtained by Construction A from a $\Z_4$-code $\calC$.\\
$\overline{\Aut(\calC)}$& The subgroup of $\Sym_n$ induced from $\Aut(\calC)$ for a $\Z_4$-code $\calC$.\\
$A.B$& A group extension with normal subgroup $A$ and quotient $B$.\\
$E(m)$& The set of all weight $m$ codewords in a binary code $E$.\\
$\calC$& A $\Z_4$-code of length $n$.\\
$\calC_0$& $\calC_0=\{ (\al_1,\dots ,\al_n)\in \{0,1\}^n|\, (2\al_1, \dots,
2\al_n)\in \calC\}$.\\
$\calC_1$& $\calC _1=\{(\al_1,\dots,\al_n)\mod2 |\ (\al_1,\dots,\al_n)\in\calC\}$.\\

$C$& The binary code of length $r$ defined by the $T_r$-module structure of $V^0$.\\
\end{tabular}
\end{center}
\begin{center}
\begin{tabular}{ll}
$D$& The binary code of length $r$ consisting of $\beta\in\Z_2^r$ with $V^\beta\neq0$.\\
$d$& The linear map from $\Z_2^n$ to $\Z_2^{2n}$, $(c_1,c_2,\dots,c_n)\mapsto(c_1,c_1,c_2,c_2,\dots,c_n,c_n)$.\\
$e$& The linear map from $\Z_2^n$ to $\Z_2^{2n}$, $(c_1,c_2,\dots,c_n)\mapsto(0,c_1,0,c_2,0,\dots,0,c_n)$.\\
$E_8$& The root lattice of type $E_8$.\\
$\mathcal{E}_n$& The binary code of length $n$ consisting of all even weight codewords.\\
$F$& A subset $\{\alpha_1,\dots,\alpha_n\}$ of a lattice $L$ of rank $n$ such that $(\alpha_i,\alpha_j)=4\delta_{ij}$, \\ & or the sublattice $\oplus_{i=1}^n\Z\alpha_i$.\\
$K$& The quotient $\stab_{\Aut(V)}(T_{2n})/\pstab_{\Aut(V)}(T_{2n})$\\ & for $V=V_L$ or $V=\tilde{V}_L$ associated to an even lattice $L$ of rank $n$.\\
$\Lambda$& The Leech lattice.\\

$\Omega_n$& An $n$-set $\{1,2,\dots,n\}$.\\
$\varphi_k$& The natural epimorphism from $\Z^n$ to $\Z_k^n$,\\ & or that from $\Z_l^n$ to $\Z_k^n$ if $k$ divides $l$.\\
$T_r$& A subVOA of a VOA with rank $r/2$ isomorphic to the tensor product of\\ & $r$-copies of $L(\shf,0)$, or a set of $r$ mutually orthogonal Ising vectors.\\
$\stab_{\Aut(V)}(T_r)$& The subgroup of $\Aut(V)$ which fixes the Virasoro frame $T_r$ setwise.\\
$\pstab_{\Aut(V)}(T_r)$& The subgroup of $\Aut(V)$ which fixes the Virasoro frame $T_r$ pointwise.\\
$\Sym_n$& The symmetric group of degree $n$.\\
$V^\beta$& The sum of irreducible $T_r$-submodules of $V$ isomorphic to \\&
$\otimes_{i=1}^r L(\shf,h_i)$ with $h_i=1/16$ if and only if $\beta_i=1$.\\
$V_L$& The lattice VOA associated to an even lattice $L$.\\
$V_L^+$& The subVOA of $V_L$ consisting of vectors fixed by the lift of $-1\in\Aut(L)$.\\
$V_E$& The code VOA associated to a binary code $E$.\\
$\tilde{V}_L$& The VOA obtained by the $\Z_2$-orbifold construction from $V_L$\\ & associated to an even unimodular lattice $L$.\\
$V^\natural$& The moonshine VOA.\\
$|X|$& The (Hamming) weight of an element $X$ of $\Z_2^n$.\\
$\Z_k$& The set of integers modulo $k$.\\
$\Z_2^n$& An $n$-dimensional vector space over $\Z_2$, or the power set of $\Omega_n$.\\

\end{tabular}
\end{center}

\section{Framed vertex operator algebra}
In this section, we review some basic facts about framed VOAs from \cite{DGH,M3}.
Every vertex operator algebra is defined over the complex number field $\C$ unless otherwise stated.
For the detail of VOAs, see \cite{Bo,FLM,FHL}.

\begin{definition}\label{df:3.1}
  A Virasoro vector $e$ is called an {\it Ising vector} if the subalgebra $\vir(e)$ generated by $e$ is isomorphic to the simple Virasoro VOA $L(\shf,0)$.
  Two Virasoro vectors $u,v\in V$ are said to be {\it orthogonal} if
  $[Y(u,z_1), Y(v,z_2)]=0$.
  A decomposition $\w=e^1+\cds +e^r$ of the conformal vector $\w$ of $V$ is
  said to be {\it orthogonal} if  $e^i$ are mutually orthogonal Virasoro vectors.
\end{definition}

\begin{remark} It is well-known that $L(\shf,0)$ is rational,
$C_2$-cofinite and has three irreducible modules $L(\shf,0)$,
$L(\shf,\shf)$ and $L(\shf,\sfr{1}{16})$. The fusion rules of
$L(\shf,0)$-modules are computed in \cite{DMZ}:
\begin{equation}\label{eq:3.1}
\begin{array}{l}
  L(\shf,\shf)\fusion L(\shf,\shf)=L(\shf,0),
  \q
  L(\shf,\shf)\fusion L(\shf,\sfr{1}{16})=L(\shf,\sfr{1}{16}),
  \vsb\\
  L(\shf,\sfr{1}{16})\fusion L(\shf,\sfr{1}{16})=L(\shf,0)\oplus L(\shf,\shf).
\end{array}
\end{equation}
\end{remark}

\begin{definition}\label{df:3.2}
  (\cite{DGH})
  A simple VOA $(V,\w)$ is said to be \emph{framed}
  if there exists a set $\{e^1, \dots,e^r\}$ of Ising vectors of $V$ such that
  $\w=e^1+\cds +e^r$ is an orthogonal decomposition.
  The full subVOA  $T_r$ generated by $e^1,\dots,e^r$ is called
  an \emph{Virasoro frame} or simply  a \emph{frame} of $V$.
  By abuse of notation, we sometimes call the set of Ising vectors
  $\{ e^1,\dots, e^r\}$ a {\it frame}, also.
\end{definition}

Given a framed VOA $V$ with a frame $T_r$, one can associate two binary
codes $C$ and $D$ of length $r$ to $V$ and $T_{r}$ as follows:
Since $T_r= L(\shf,0)^{\otimes r}$ is rational, $V$ is a completely reducible $T_r$-module. That is,
\begin{equation*}\label{2.1}
V \cong \oplus_{h_i\in\{0,\frac{1}{2},\frac{1}{16}\}}
m_{h_1,\ldots, h_{r}}L(h_1,\ldots,h_{r}),
\end{equation*}
where $L(h_1,h_2,\dots,h_r)=L(\shf,h_1)\otimes\cdots\otimes L(\shf, h_r)$ and the nonnegative integer $m_{h_1,\ldots,h_r}$ is the
multiplicity of $L(h_1,\ldots,h_r)$ in $V$. In particular, all the
multiplicities are finite and $m_{h_1,\ldots,h_r}$ is at most $1$ if
all $h_i$ are different from $\sfr{1}{16}$.

Let $U\cong L(h_1,h_2,\dots,h_r)$ be an irreducible
module for $T_{r}$. The $\tau$-word $\tau(U)$ of $U$ is a binary word
$\beta=(\be_1, \dots,\be_r)\in \Z_2^r$ such that
\begin{equation}
\be_i=
\begin{cases}
0 & \text{ if } h_i=0 \text{ or } \shf,\\
1 & \text{ if } h_i=\sfr{1}{16}.
\end{cases}
\end{equation}
For any $\beta\in \Z_2^r$, define $V^\be$ as the sum of all
irreducible submodules $U$ of $V$ such that $\tau(U)=\be$. Set $$
D:=\{\be\in \Z_2^r \mid V^\be \ne 0\}. $$ Then $D$ is an even linear
code of length $r$ and $V=\oplus_{\be\in D} V^\be$.

For any $c=(c_1,...,c_r)\in \Z_2^r$, denote
$M^c=m_{h_1,...,h_r}L(h_1,...,h_r)$  where $h_i=\shf$ if
$c_i=1$ and $h_i=0$ elsewhere. Set
 \begin{equation*}\label{2.3}
C:=\{c\in \Z_2^r \mid M^c \neq 0\}.
\end{equation*}
Then $V^0=\oplus_{c\in C}M^c$ is the code VOA $V_C$ associated to $C$ (\cite{M4}).

Summarizing, there exists a pair $(C,D)$ of even linear codes such
that $V$ is a $D$-graded extension of a code VOA $V_C$ associated
to $C$. The  pair of codes $(C,D)$ is called the {\it structure codes} of a
framed VOA $V$ associated to the frame $T_{r}$. Since the powers of
$z$ in an $L(\shf,0)$-intertwining operator of type
$L(\shf,\shf)\times L(\shf,\sfr{1}{16})\to L(\shf,\sfr{1}{16})$ are
half-integral, the structure codes $(C, D)$ satisfy $C\subset
D^\perp$.  Moreover, the following theorem holds (cf. \cite[Theorem 2.9]{DGH} and \cite[Theorem 6.1]{M3}).
\begin{theorem}
Let $V$ be a framed VOA with structure codes $(C,D)$. Then, $V$ is
holomorphic if and only if $C=D^\perp$.
\end{theorem}

\begin{remark}\label{rem:2.5}
Let $V$ be a framed VOA with structure codes $(C,D)$, where
$C,D\subset \Z_2^r$. For a binary codeword $\beta \in \Z_2^r$, we
define
\begin{equation}\label{eq:3.3}
  \ds \tau_\beta(u):=(-1)^{\la \al,\be\ra} u \q \text{ for } u\in
  V^\al.
\end{equation}
Then by the fusion rules, $\tau_\be$ defines an automorphism on $V$
\cite[Theorem 4.7]{M1}. Note that the subgroup $\mathcal{T}=\{ \tau_\be\,|\, \be\in
\Z_2^r\}$ is an elementary abelian $2$-group and is isomorphic to
$\Z_2^r/D^\perp$. In addition, the fixed point subspace $V^{\mathcal{T}}$ is
equal to $V^0$ and all $V^\al,\al\in D$ are irreducible
$V^0$-modules. Similarly, we can define an automorphism on $V^0$ by
\[
\sigma_\beta(u):=(-1)^{\la \al,\be\ra} u \q \text{ for } u\in
  M^\al,
\]
where $V^0=\oplus_{\al\in C}M^\al$. Note that the group $Q=\{
\sigma_\be\,|\, \be\in \Z_2^r\}\cong \Z_2^r/C^\perp$ is elementary
abelian and $(V^0)^Q=M^{0}=T_r$.
 \end{remark}

\section{Frame stabilizers}

In this section, we shall recall the definitions of frame stabilizers and pointwise frame stabilizers of a framed VOA. Some basic properties will also be reviewed from \cite{DGH,LY}.

\begin{definition}\label{df:6.1}
  Let $V$ be a framed VOA with an Virasoro frame $T_r=\vir(e^1)\tensor \cds \tensor \vir(e^r)$.
  The {\it frame stabilizer} of $T_{r}$ is the subgroup of $\aut(V)$ which stabilizes
  the frame $T_{r}$ setwise. The {\it pointwise frame stabilizer} is
  the subgroup of $\aut(V)$ which fixes $T_{r}$ pointwise.
  The frame stabilizer and the pointwise frame stabilizer of $T_{r}$ are denoted by $\stab_{\aut(V)} (T_{r})$ and
  $\pstab_{\aut(V)}(T_{r})$, respectively.
\end{definition}

Let $(C,D)$ be the structure codes of $V$ with respect to $T_{r}$, i.e.,
$$
  V=\oplus_{\alpha\in D} V^\alpha,\quad  \tau(V^\alpha)=\alpha \quad
  \text{ and}\quad V^0= V_C.
$$

If $D=0$, i.e., $V=V_C$, then the frame stabilizer can be determined easily.

\begin{lemma}
Let $C$ be an even linear code and $Q=\{\sigma_\gamma | \, \gamma\in \Z_2^r \}\subset \aut(V_C)$. Then we have an exact sequence
\[
1 \longrightarrow Q \longrightarrow \stab_{\aut(V_C)}(T_{r}) \bar\longrightarrow \aut(C) \longrightarrow 1.
\]
\end{lemma}
\begin{proof} Let $g\in \stab_{\aut(V_C)}(T_{r})$. Then $g$ gives a permutation $\bar{g}$ on the set of Ising vectors $\{e_1,e_2,\dots, e_r\}$.
It is easy to see that for any $c\in C$, $g(M^c)=M^{\bar{g}(c)}$.
Hence $\bar{g}\in\aut(C)$.
By \cite[Lemma 3.4]{M1}, the map $\ \bar{}\ $ is surjective.
If $\bar{g}=1$ then $g$ is trivial on $T_{r}$.
Hence the kernel of $\ \bar{}\ $ is equal to the pointwise stabilizer $\pstab_{\aut(V_C)}(T_{r})$ of $T_{r}$, and which is equal to $Q$ (see \cite[Section 6]{LY}).
\end{proof}

\medskip

Next, let us review the properties of the pointwise frame stabilizers for framed VOAs from \cite[Section 6]{LY}.

\begin{notation}
For $\alpha=(\alpha_1,\dots,\alpha_r), \beta\in
(\beta_1,\dots,\beta_r)\in \Z_2^r$, we define
$$
  \alpha\cd \beta=\alpha\cap\beta:= (\alpha_1\beta_1,\dots,\alpha_r\beta_r) \in \Z_2^r .
$$
That is, the product $\alpha\cd \beta$ is taken in the ring $\Z_2^r$.
\end{notation}

\begin{theorem}\cite[Theorem 12]{LY}\label{thm:6.7}
  Let $V$ be a framed VOA with structure codes  $(C,D)$. Let  $\xi\in \Z_2^r \setminus C^\perp$. Then, there exists
  $\theta\in \pstab_{\Aut(V)}(T_r)$ such that $\theta|_{V^0}=\sigma_{\xi}$
  if and only if $\alpha\cd \xi\in C$ for all $\alpha \in D$. Moreover, $\theta$ has order $2$ if $\wt(\alpha\cd \xi) \equiv 0 \mod 4$ for all $\alpha \in D$; otherwise, $\theta$ has order $4$.
\end{theorem}

Define $P:=\{ \xi \in \Z_2^r\mid \alpha\cd \xi \in C\ \text{for all}\ \alpha\in D\}$.
It is clear that $P$ is a linear subcode of $C$. Moreover, we have the following.

\begin{theorem}\cite[Theorem 13]{LY}\label{thn:5.5}
Let $V$ be a framed VOA with structure codes  $(C,D)$ and let $\mathcal{T} =\{\tau_\alpha|\ \alpha\in \Z_2^r\}$ be the subgroup generated by the Miyamoto involutions associated to Ising vectors in the frame $T_r$. Then we have the following central extension:
$$
\begin{array}{ccccccccc}
  1 & \longrightarrow &\mathcal{T} & \longrightarrow & \pstab_{\Aut(V)}(T_r)
  & \longrightarrow & \pstab_{\Aut(V)}(T_r)/\mathcal{T} & \longrightarrow &  1
  \vsb\\
   && \downarrow\! \wr && ||  && \downarrow\! \wr &&
  \vsb\\
  1 &\longrightarrow& \Z_2^r/D^\perp &\longrightarrow& \pstab_{\Aut(V)}(T_r)
  &\longrightarrow&  P/C^\perp &\longrightarrow& 1
\end{array}
$$
The commutator relation in $\pstab_{\Aut(V)}(T_r)$ can  be described as follows.

For $\xi^1,\xi^2\in P$, let $\theta_{\xi^i}, i=1,2,$ be an extension of $\sigma_{\xi^i}$
  to $\pstab_{\Aut(V)}(T_r)$.
  Then $[\theta_{\xi^1},\theta_{\xi^2}]=1$ if and only if
  $\la \alpha\cd \xi^1,\alpha\cd \xi^2\ra= 0$ for all $\alpha \in D$.
\end{theorem}

The full frame stabilizer is much more complicated. Nevertheless, we still have the following.
\begin{lemma}{\rm (cf. \cite[Theorem 2.8]{DGH} and \cite[Section 6]{LY})}\label{Lem 5.6} Let $V$ be a framed VOA with structure codes $(C,D)$.
Then the quotient group $\stab_{\aut(V)}(T_r)/\pstab_{\aut(V)}(T_r)$ is isomorphic to a subgroup of $\aut(C)\cap\aut(D)$.
\end{lemma}

\section{$\Z_4$-codes, lattices and framed VOAs}
In this section, we review $\Z_4$-codes and lattices obtained by Construction A from $\Z_4$-codes from \cite{CS98}.
Moreover, we review the structure codes of framed VOAs associated to even lattices having $4$-frames from \cite{DGH}.

A subgroup $\calC$ of $\Z_4^n$ is called a (linear) {\it $\Z_4$-code} of length $n$.
The {\it (Euclidean) weight} of $c=(c_1,c_2,\dots,c_n)\in\Z_4^n$, $c_i\in\{0,\pm1,2\}$, is $\sum_{i=1}^n c_i^2$.
The {\it dual code} of $\calC$ is defined as $\calC^\perp=\{x\in\Z_4^n\mid \langle x,y\rangle=0,\text{ for all } y\in\calC\}$, where $\langle x,y\rangle=\sum_{i=1}^nx_iy_i\in\Z_4$.
A $\Z_4$-code $\calC$ is said to be {\it self-orthogonal} if $\calC\subset\calC^\perp$, and is said to be {\it self-dual} if $\calC=\calC^\perp$.
A self-dual $\Z_4$-code $\calC$ is said to be {\it Type II} if the Euclidean weight of any element in $\calC$ is divisible by $8$.
A Type II $\Z_4$-code $\calC$ is said to be {\it extremal} if the minimum weight of $\calC$ is equal to $8(\lfloor n/24\rfloor+1)$.
The automorphism group of $\calC$ is the subgroup of $\aut(\Z_4^n)\cong R:\Sym_n$ preserving $\calC$, where $R$ consists of sign change maps on coordinates.
Hence we obtain the following exact sequence
\begin{eqnarray}
1\to (R\cap\aut(\calC))\to\aut(\calC)\to\overline{\aut(\calC)}\to1,\label{Def:C^0}
\end{eqnarray}
where $\overline{\aut(\calC)}\cong\aut(\calC)/(R\cap\aut(\calC))$ is a subgroup of $\Sym_n$.

\medskip

Now let us study the structure codes for the lattice VOA $V_{L}$ associated to $L$ having $4$-frame.

\begin{definition}
Let $L$ be an even lattice of rank $n$. A subset $\{
\al_1, \dots, \al_n\}$ is called \textit{a $4$-frame} of $L$ if $(\al_i, \al_j)=4\delta_{i,j}$ for all $i,j$.
By abuse of notation, we sometimes call the sublattice $\oplus_{i=1}^n\Z \al_i$ a $4$-frame, also.
\end{definition}

Let $L$ be an even lattice of rank $n$ and $V_L$ the VOA associated to $L$ (\cite{Bo,FLM}).
Let $\al\in L$ with $(\al,\al)=4$. It is well-known (cf. \cite{DMZ}) that
\[
\om^{\pm} (\al)= \frac{1}{16}\al(-1)^2\cdot 1 \pm \frac{1}4 (e^\al +e^{-\al} )
\]
are two mutually orthogonal Ising vectors in $V_L$.
If $L$ contains a $4$-frame $F=\{\al_1, \dots, \al_n\}$, then the lattice VOA $V_L$ is a framed VOA with a Virasoro frame $$\{\om^+(\al_1), \om^-(\al_1) , \dots, \om^+(\al_n), \om^-(\al_n)\}.$$
We call it {\it the Virasoro frame associated to} $F$.

\begin{notation}
For any positive integer $k$, denote the natural epimorphism from $\Z$ to $\Z_k$ by $\varphi_k$. We shall also extend $ \varphi_k$ to a homomorphism from $\Z^n$ to $\Z_k^n$ by
\[
\varphi_k(a_1, \dots, a_n)= (\varphi_k a_1, \dots, \varphi_k a_n).
\]
By abuse of notation, we also use $\varphi_k$ to denote the natural projection
from $\Z_\ell^n$ to $\Z_k^n$ when $k$ divides $\ell$.
\end{notation}

Let $\{\al_1, \dots, \al_n\}$ be a $4$-frame of a lattice $L$ and $F=\Z\al_1\oplus \cdots\oplus \Z\al_n$. Then $L/F \subset F^*/F \cong \Z_4^n$ forms a $\Z_4$-code.

Let $\calC$ be a $\Z_4$-code of length $n$ such that $\calC \cong L/F$. Then
\begin{eqnarray}
L\cong A_4(\calC)=\frac{1}2\left\{ (x_1,\dots,x_n)\in \mathbb{Z}^n|\,
\varphi_4(x_1,\dots,x_n)\in \calC \right\}.\label{Eq:A4}
\end{eqnarray}
The following lemma is well-known.

\begin{lemma}[cf. \cite{BSC}]
\begin{enumerate}
\item $A_4(\calC)$ is integral if and only if $\calC$ is self-orthogonal.
\item $A_4(\calC)$ is even if and only if the Euclidean weight of
any element in $\calC$ is divisible by $8$.

\item $A_4(\calC)$ is even unimodular if and only if $\calC$ is type II.
\end{enumerate}
\end{lemma}

\begin{remark}
 Note that $A_4(0)\cong\Z\al_1\oplus \cdots\oplus \Z\al_n$, $(\alpha_i,\alpha_j)=4\delta_{ij}$ and it gives a $4$-frame of $A_4(\calC)$.
Hence $\{\om^+(\al_1), \om^-(\al_1), \dots, \om^+(\al_n), \om^-(\al_n)\}$ defines a Virasoro frame for $V_{A_4(\calC)}$. We shall call it the {\it Virasoro frame associated to $\calC$}.
\end{remark}

Let $\calC $ be a $\Z_4$-code of length $n$. We shall define two binary codes
\[
\begin{split}
&\calC_0=\{ (\al_1,\dots ,\al_n)\in \{0,1\}^n|\, (2\al_1, \dots,
2\al_n)\in \calC\}, \\
& \calC_1=\{ \varphi_2(\al)\ |\ \al\in \calC \}=\{(\al_1,\dots,\al_n)\mod2 |\ (\al_1,\dots,\al_n)\in\calC\}.\\
\end{split}
\]
Note that if we define a linear map $t: \calC \to \calC$ by $t(\al)=2\al$, then $\calC_1\cong \calC/ \ker t$ and $\ker t\cong \calC_0$ as abelian groups.  Thus, as an abelian group, $\calC$ is an extension of $\calC_1$ by $\calC_0$  and we have the exact sequence
\[
0\to \calC_0 \to \calC\to \calC_1 \to 0.
\]
Note also that $\calC_1\subset\calC_0$ and if $|\calC_1|=2^{k_1}$ and $|\calC_0|=2^{k_0}$, then the original $\Z_4$-code $\calC\cong \Z_2^{k_0-k_1}\times \Z_4^{k_1}$ as an abelian group.

If $\calC$ is self-orthogonal then both $\calC _0$ and $\calC_1$ are even binary codes of length $n$, and $\calC_0\subset \calC_1^\perp$.
If $\calC$ is self-dual then $\calC_0=\calC_1^\perp$, and if the Euclidean weight of any element in $\calC$ is divisible by $8$ then $\calC_1$ is doubly even.
By the definitions of $\calC_0$ and $\calC_1$, the group $\overline{\aut(\calC)}$ preserves both $\Aut(\calC_0)$ and $\Aut(\calC_1)$, that is, $\overline{\aut(\calC)}\subset\Aut(\calC_0)\cap\Aut(\calC_1)$.

\begin{remark}
The binary codes $\calC_0$ and $\calC_1$ may also be thought of as $\Z_4$-code analogues of the structure codes $C$ and $D$. In fact, if $L$ is an even lattice with a $4$-frame $F$ and  $L/F\cong \calC$ as a $\Z_4$-code, then the structures $(C,D)$ for the lattice VOA $V_L$  (and for its $\Z_2$-orbifold $\tilde{V}_L$) with respect to the Virasoro frame associated to $\calC$ are closely related to $\calC_0$ and $\calC_1$ (see Proposition \ref{PDGH} and \ref{Prop DGH}).
\end{remark}

\begin{notation}\label{Nde}
Let $d$ and $e$ denote the linear maps from $\Z_2^n$ to $\Z_2^{2n}$ defined by
\[
d:(c_1,c_2,\dots,c_n)\mapsto(c_1,c_1,c_2,c_2,\dots,c_n,c_n)
\]
and
\[
 e:(c_1,c_2,\dots,c_n)\mapsto(0,c_1,0,c_2,0,\dots,0,c_n).
\]
\end{notation}

The structure codes of $V_{A_4(\calC)}$ with respect to the Virasoro frame associated to $\calC$ is described in \cite{DGH} as follows.

\begin{proposition}\label{PDGH} \cite[Corollary 3.3]{DGH}
Let $\calC $ be a $\Z_4$-code of length $n$ such that $A_4(\calC)$ is even, and $\calC _0$ and $\calC _1$
the binary codes defined as above. Then the structure codes of the lattice VOA
$V_{A_4(\calC )}$ with respect to the Virasoro frame associated to $\calC$ are given by
\[
D=d(\calC _1)\quad \text{ and }\quad  C=\Span_{\Z_2}\{ d(\Z_2^n), e(\calC_0)\}.
\]
\end{proposition}

\medskip

Now let $L$ be an even unimodular lattice having a $4$-frame $F$ and let $\theta\in \aut(V_L)$ be a lift of $-1\in\aut(L)$.
Let $V_L^+$ denote the subVOA of $V_L$ consisting of vectors in $V_L$ fixed by $\theta$.
By the definition of Ising vectors associated to $F$, the Virasoro frame of $V_L$ associated to $F$ is contained in $V_L^+$. Hence $V_L^+$ is framed.

Let $V_L^T$ be the unique $\theta$-twisted module for $V_L$ and $V_L^{T,+}$ the irreducible $V_L^+$-submodule of $V_L^T$ with integral weights. 
Let $$ \tilde{V}_L = V_L^+ \oplus V_L^{T,+}.$$ It is known that $ \tilde{V}_L$
has a unique VOA structure by extending its
$V_L^+$-module structure (see \cite{FLM} and \cite[Proposition 8]{LY}).
The construction of VOA  $\tilde{V}_L$ is often called a {\it $\Z_2$-orbifold construction} from $V_L$ .
Clearly, $\tilde{V}_L$ is also framed.

\begin{notation}\label{NEn}
 Denoted by $\mathcal{E}_n$  the code consisting of all codewords in $\Z_2^n$ of even weight. Then the code  $d(\mathcal{E}_n)$ is generated by
\[
\begin{pmatrix}
1111 &0000& 0000& \cdots &0000\\
0011 &1100& 0000& \cdots &0000\\
0000 &1111& 0000& \cdots &0000\\
0000 &0011& 1100& \cdots &0000\\
\vdots & & \vdots &\ddots & \vdots\\
0000 & 0000& 0000 & \cdots&1111
\end{pmatrix}_.
\]
\end{notation}

\begin{proposition}\cite[Corollaries 3.3 and 3.5]{DGH}\label{Prop DGH}
Let $\calC $ be a type II $\Z_4$-code of length $n$, and $\calC _0$ and $\calC _1$
the associated binary codes.
Then the lattice $A_4(\calC)$ is even unimodular and the structure codes of $\tilde{V}_{A_4(\calC)}$ with respect to the Virasoro frame associated to $\calC$ are given by
\[
D= \Span_{\Z_2}\{ d(\calC _1), e((1^n))\}\quad \text{ and }\quad  C=\Span_{\Z_2}\{ d(\mathcal{E}_n) , e(\calC_0)\}.
\]
\end{proposition}

\section{Frame stabilizer of the frame of $V_L$ associated to a $4$-frame of $L$}
Let $L$ be an even lattice of rank $n$ having $4$-frame $F$.
Let $T_{2n}$ be the Virasoro frame of $V_L$ associated to $F$.
In this section, we study the subgroup $K= \stab_{\aut(V_L)}(T_{2n})/\pstab_{\aut(V_L)}(T_{2n})$ of $\aut(C)$.
For the detail of the structure codes $(C,D)$ of $V_L$ associated to $T_{2n}$, see Proposition \ref{PDGH}.

\subsection{Method of calculating the frame stabilizer}
In this subsection, we discuss the subgroup $K$ of $\aut(C)$.
First, we recall triality automorphisms of $V_{A_1\oplus A_1}$ from \cite{FLM}.

\begin{lemma}\label{L0} Let $L=\Z\alpha_1\perp\Z\alpha_2$ be an even lattice of rank $2$ such that $(\alpha_1,\alpha_1)=(\alpha_2,\alpha_2)=2$ and $(\alpha_1,\alpha_2)=0$.
Let $T_{4}$ be the Virasoro frame in $V_L$ associated to the $4$-frame $F=\{\alpha_1\pm\alpha_2\}$.
Then the stabilizer $\stab_{N(V_L)}(T_4)$ of $T_{4}$ in $N(V_L)$ acts on $T_4$ as $\Sym_4$, where $N(V_L)$ is the normal subgroup of $\aut(V_L)$ generated by $\{\exp(v_0)\mid v\in (V_L)_1\}$.
\end{lemma}
\begin{proof} Set $e_1=\omega^+(\alpha_1+\alpha_2)$, $e_2=\omega^-(\alpha_1+\alpha_2)$, $e_3=\omega^+(\alpha_1-\alpha_2)$, $e_4=\omega^-(\alpha_1-\alpha_2)$.
Then $T_4=\{e_1,e_2,e_3,e_4\}$, and the structure codes of $V_L$ associated to $T_4$ are $(\mathcal{E}_4,0)$.
By Lemma \ref{Lem 5.6}, there is a canonical group homomorphism $\psi:\stab_{N(V_L)}(T_4)\to \aut(\mathcal{E}_4)\cong\Sym_4$.
We now identify $e_i$ with $i$.
Then $\exp(\frac{\pi i}{4}(\alpha_1\pm\alpha_2)(0))\in N(V_L)$ act as $(1\ 2)$ and $(3\ 4)$ on $T_4$.
Hence $\Imm\ \psi$ contains $\langle (1\ 2),(3\ 4)\rangle$.

By \cite[Corollary 11.2.4]{FLM}, there is an automorphism $\sigma\in N(V_L)$ such that
\begin{eqnarray*} \sigma(\alpha_i(-1)^2)&=&\alpha_i(-1)^2,\\
\sigma(\alpha_1(-1)\alpha_2(-1))&=&e^{\alpha_1+\alpha_2}+\theta(e^{\alpha_1+\alpha_2})+e^{\alpha_1-\alpha_2}+\theta(e^{\alpha_1-\alpha_2}),\\
\sigma(e^{\alpha_1+\alpha_2}+\theta(e^{\alpha_1+\alpha_2})+e^{\alpha_1-\alpha_2}+\theta(e^{\alpha_1-\alpha_2}))&=&\alpha_1(-1)\alpha_2(-1),\\
\sigma (e^{\alpha_1+\alpha_2}+\theta (e^{\alpha_1+\alpha_2})-(e^{\alpha_1-\alpha_2}+\theta (e^{\alpha_1-\alpha_2})))&=&e^{\alpha_1+\alpha_2}+\theta(e^{\alpha_1+\alpha_2})-(e^{\alpha_1-\alpha_2}+\theta (e^{\alpha_1-\alpha_2})).
\end{eqnarray*}
Then one can see that
\begin{eqnarray*}
\sigma(e_1)=e_1, \sigma(e_2)=e_3, \sigma(e_3)=e_2,\ \sigma(e_4)=e_4.
\end{eqnarray*}
Hence $\Imm\ \psi$ contains $(2\ 3)$.
Thus $\Imm\ \psi\cong\Sym_4$.
\end{proof}

\begin{lemma}\label{L1} Let $L$ be an even lattice of rank $n$ having $4$-frame $F$.
Let $\calC=L/F$ be a $\Z_4$-code of length $n$.
Let $T_{2n}$ be the Virasoro frame of $V_L$ associated to $F$. Then the stabilizer in $K$ of the subcode $d(\Z_2^n)$ of $C$ is isomorphic to $2\wr\overline{\aut(\calC)}$.
\end{lemma}

\begin{proof}
Let $S$ denote the stabilizer of $d(\Z_2^n)$ in $K$.

For $k=1, \dots,n$, define $t_k= \exp(\frac{\pi i \al_k(0)}{4})$.
Then $t_k$ maps $\om^{\pm}(\al_k) $ to $\om^{\mp}(\al_k)$ and fixes $\om^\pm(\al_i)$ if $i\neq k$. Thus, $ t_k\in \stab_{\aut(V_L)}(T_{2n})$ and it stabilizes $d(\Z_2^n)$. Let $\bar{t}_k$ be the image of $t_k$ under the canonical homomorphism $\stab_{\aut(V_L)}(T_{2n})\to K$. Then $\bar{t}_1, \dots, \bar{t}_n$ generate a subgroup $A\cong 2^{n}$ in $S$.
Let $g\in\Aut(\calC)$.
By (\ref{Eq:A4}), $g$ induces an automorphism of $L$ preserving $F$.
Hence by \cite[Proposition 5.4.1]{FLM} it lifts to an automorphism $g^\prime$ of $V_L$ preserving $T_{2n}$.
Since $R\cap\Aut(\calC)$ is trivial on $T_{2n}$, $g^\prime$ acts as $g(R\cap\Aut(\calC))\in\overline{\aut(\calC)}$ on $\{\{\omega^\pm(\alpha_i)\}\mid i\in\Omega_n\}$.
Thus we obtain a subgroup $A.\overline{\aut(\calC)}$ of $S$.

Now suppose $g\in \stab_{\Aut(V_L)}(T_{2n})$ such that $g\cdot \pstab_{\Aut(V_L)}(T_{2n})\in S$.
Since the structure codes of $V_F$ associated to $T_{2n}$ are $(d(\Z_2^n),0)$, $g$ stabilizes $V_F$.
Since $F$ has no roots, we have $\Aut(V_F)\cong (\C^\times)^n.\Aut(F)$ (cf.\ \cite{DN}).
We view $g$ as an element in $\stab_{\aut(V_F)}(T_{2n})$.
Note that $(\C^\times)^n\cap \stab_{\aut(V_F)}(T_{2n})=A$.
Set $\bar{g}=gA\in\aut(F)$.
It follows from $g\in\aut(V_L)$ that $\bar{g}\in\aut(L)$.
Hence $g$ induces an automorphism of $\calC\cong L/F$, modulo the action of $\pstab_{\aut(V_L)}(T_{2n}). A$, and $S$ is a subgroup of $A.\overline{\aut(\calC)}\cong 2\wr \overline{\aut(\calC)}$.
Therefore $S\cong 2\wr \overline{\aut(\calC)}$.
 \end{proof}

\begin{remark} Lemma \ref{L1} was essentially proved by \cite[Theorem 2.8 (iii)]{GH} in terms of lattices.
\end{remark}

\begin{proposition}\label{P1} Let $L$ be an even lattice of rank $n$ having $4$-frame $F$.
Let $T_{2n}$ be the Virasoro frame of $V_L$ associated to $F$ and $(C,D)$ the structure codes associated to $T_{2n}$.
Let $c=\{i,j\}$ be an element in $C(2)$.
Then there exists an element $g$ in $\stab_{\aut(V_L)}(T_{2n})$ acting as $g=(i\ j)$ on $T_{2n}$.
\end{proposition}
\begin{proof}
If $c$ belongs to $d(\Z_2^n)$, then it is clear from Lemma \ref{L1}.
Assume that $c\notin d(\Z_2^n)$.
Then there exist exactly two codewords $c_1=\{i,k\},c_2=\{j,l\}\in d(\Z_2^n)$ such that $|c\cap c_1|=|c\cap c_2|=1$ and $k\neq l$.
We view $Z=\Span_{\Z_2}\{ c_1,c_2,c\}$ as a binary code of length $4$ on $\{i,j,k,l\}$.
Then $Z\cong\mathcal{E}_4$, and we can find $\alpha_1,\alpha_2\in L$ of norm $2$ such that $V_Z=V_{\Z\alpha_1\oplus\Z\alpha_2}\cong V_{A_1\oplus A_1}$ and $\alpha_1\pm\alpha_2\in F$.
By Lemma \ref{L0}, the canonical group homomorphism $\stab_{N(V_Z)}(\{e_i,e_j,e_l,e_k\})\to \Sym(\{i,j,k,l\})$ is surjective.
Note that $\stab_{N(V_Z)}(\{e_i,e_j,e_l,e_k\})$ lifts to a subgroup of $N(V_L)$ preserving $T_{2n}$.
Hence there exists an automorphism of $V$ acting as $(i\ j)$ on $T_{2n}$.
\end{proof}

For an element $h=\{i,j\}\in C(2)$, let $\xi(h)$ denote the element of $K$ acting as $(i\ j)\in\Sym(T_{2n})\cong \Sym_{2n}$.
\begin{lemma}\label{L2} Let $H$ and $H^\prime$ be subcodes of $C$ isomorphic to $d(\Z_2^n)$.
Then there exists $g\in \langle\xi(h)\mid h\in \Span_{\Z_2}\{ H, H^\prime\}(2)\rangle\subset K$ such that $g(H)=H^\prime$.
\end{lemma}
\begin{proof} Assume $H\neq H^\prime$.
Let $h\in H(2)\setminus H^\prime$.
Then there is $h^\prime\in H^\prime(2)$ such that $|h\cap h^\prime|=1$.
Clearly, $h^\prime\notin H(2)$.
It is easy to see that $\xi(h)\circ \xi(h^\prime)(h)=h^\prime$. Set $g=\xi(h)\circ \xi(h^\prime)$.
Since $H\cong d(\Z_2^n)$, $H(2)\setminus \{h\}$ is fixed by $\xi(h)$ pointwisely. Similarly, $H'(2)\setminus \{h'\}$ is fixed by $\xi(h')$ pointwisely. Thus, $g(H(2)\cap H^\prime (2))=H(2)\cap H^\prime (2)$. On the other hand, $h^\prime\in g(H)(2)\cap H^\prime (2)$ and hence $|g(H)(2)\cap H^\prime(2)|>|H(2)\cap H^\prime(2)|$.
Note that $g(H)\subset C$ as $g\in \aut(C)$.
By reverse induction on $|H(2)\cap H^\prime(2)|$, we are done.
\end{proof}

Let $\mathcal{H}$ denote the set of all subcodes of $C$ isomorphic to $d(\Z_2^n)$.
By Lemmas \ref{L1}, \ref{L2} and Proposition \ref{P1}, we obtain the following theorem.

\begin{theorem}\label{Thm7.6} Let $L$ be an even lattice of rank $n$ having $4$-frame $F$.
Let $T_{2n}$ be the Virasoro frame of $V_L$ associated to $F$ and $(C,D)$ the structure codes associated to $T_{2n}$.
Then $K=\stab_{\aut(V_L)}(T_{2n})/\pstab_{\aut(V_L)}(T_{2n})$ is transitive on $\mathcal{H}$, and it is generated by the subgroup of shape $2\wr\overline{\aut(\calC)}$ and $\{\xi(h)\mid h\in C(2)\}$.
Moreover, $|K:2\wr\overline{\aut(\calC)}|=|\mathcal{H}|$.
\end{theorem}

Now, we consider the action of $\Aut(C)$ on $\mathcal{H}$.

\begin{lemma} \label{H} The stabilizer in $\Aut(C)$ of $d(\Z_2^n)$ is isomorphic to $2 \wr \aut(\calC_0)$.
In particular
\[
 |\mathcal{H}|=|\aut(C): 2 \wr \aut(\calC_0)|.
\]
\end{lemma}

\begin{proof}
 By Theorem \ref{Thm7.6}, $\aut( C)$ acts transitively on $\mathcal{H}$. Thus
\[
 |\mathcal{H}|=|\aut(C): \stab_{\aut(C)}(d(\Z_2^n))|.
\]
Since $C=\Span_{\Z_2}\{ d(\Z_2^n), e(\calC_0)\}$, it is easy to see that $\stab_{\aut(C)}(d(\Z_2^n))$ is isomorphic to
$2\wr \aut(\calC_0)$ and we have the desired conclusion.
\end{proof}

As a corollary, we have the following.

\begin{corollary}\label{Prop7.8}
$
 |\aut(C): K| = |\aut(\calC_0): \overline{\aut(\calC)}|.
$
\end{corollary}

\begin{remark} Let $g\in\Aut(C)$.
Then there exists $h\in K$ such that $hg$ stabilizes $d(\Z_2^n)$ by Theorem \ref{Thm7.6}.
Hence $g\in K$ if and only if $hg$ belongs to the stabilizer of $d(\Z_2^n)$ in $K$ determined in Lemma \ref{L1}.
\end{remark}

\medskip

\subsection{Examples}\label{SsEx} Next we shall discuss several examples in detail. Recall that there are exactly four type II $\Z_4$-codes of length $8$ (\cite{CS}). They are of the shape $4\cdot 2^6$, $4^2\cdot 2^4$, $4^3\cdot 2^2$, $4^4$ as groups, and their generator matrices are given as follows:
\[
\begin{pmatrix}
1111\, 1111\\
2200\, 0000\\
0220\, 0000\\
0022\, 0000\\
0002\, 2000\\
0000\, 2200\\
0000\, 0220
\end{pmatrix},
\qquad
\begin{pmatrix}
3111\, 3111\\
1111\, 2000\\
2200\, 0000\\
0220\, 0000\\
0000\, 2200\\
0000\, 0220
\end{pmatrix},
\qquad
\begin{pmatrix}
3111\, 3111\\
1111\, 2000\\
1320\, 1100\\
2200\, 0000\\
0220\, 0220
\end{pmatrix},
\qquad
\begin{pmatrix}
3111\, 3111\\
1111\, 2000\\
1320\, 1100\\
1010\, 1032
\end{pmatrix}.
\]

Let $\calC$ be a type II code of length $8$ over $\Z_4$.
Then $$A_4(\calC)=\frac{1}2\left\{ (x_1,\dots,x_8)\in \mathbb{Z}^8|\,
\varphi_4(x_1,\dots,x_8)\in \calC \right\}\cong E_8.$$
Let us determine the frame stabilizer in $\Aut(V_{E_8})$ of the Virasoro frame $T_{16}$ associated to $\calC$ for each case.

\medskip

\noindent \textbf{Case 1. ($\calC\cong 4\cdot 2^6$)}
In this case, $\calC_0= \mathcal{E}_8$, $\calC_1= \{ (1^8), (0^8)\}$ and $\aut(\calC_0)\cong \Sym_8$. The structure codes for $V_{A_4(\calC)}\cong V_{E_8}$ are then given by
\[
C={\mathcal{E}}_{16}, \qquad D= \{(1^{16}), (0^{16})\}.
\]
Then  $\aut(C)\cong \aut(D) =\Sym_{16}$ and the code $P=\{\al\in C|\, \al\cdot\be\in C \text{ for all } \be\in D\} = C$. Note that $\dim P=15$ and by Theorem \ref{thm:6.7},  we have $\pstab_{\aut(V_{E_8})}(T_{16}) \cong 2^{1+14}$.

Since $\aut(\calC)\cong 2^7:\Sym_8$, we have $\overline{\aut(\calC)}= \Sym_8$.
Hence $\overline{\aut(\calC)}\cong \aut(\calC_0).$ By Corollary \ref{Prop7.8}, $\stab_{\aut(V_{E_8})}(T_{16})/\pstab_{\aut(V_{E_8})}(T_{16}) \cong \aut(C)\cong \Sym_{16}$ and  $\stab_{\aut(V_{E_8})}(T_{16})$ has the shape $2^{1+14}.\Sym_{16}$.

\medskip

\noindent \textbf{Case 2. ($\calC\cong  4^2\cdot 2^4$)}
In this case, $\calC_0= \mathcal{E}_4\oplus \mathcal{E}_4$, $\calC_1= \Span_{\Z_2}\{ (1111\,0000), (0000\,1111)\}$ and $\aut(\calC_0)\cong \Sym_4\wr 2$.  Then the structure codes for $V_{A_4(\calC)}\cong V_{E_8}$ are given by
\[
C=\mathcal{E}_{8}\oplus \mathcal{E}_{8}, \qquad D= \Span_{\Z_2}\{ (1^8 0^8), (0^8 1^8)\}.
\]
Then $\aut(C)\cong \aut(D)= \Sym_8\wr 2$ and the code $P=\{\al\in C|\, \al\cdot\be\in C \text{ for all } \be\in D\} = C$. By Theorem \ref{thm:6.7}, $\pstab_{\aut(V_{E_8})}(T_{16}) \cong 2^{1+6}\times 2^{1+6}$.

Since $\aut(\calC)\cong 2^6.(\Sym_4\wr 2)$, we have $\overline{\aut(\calC)}= \Sym_4\wr 2$. Again, we have  $\overline{\aut(\calC)}\cong \aut(\calC_0)$. By Corollary \ref{Prop7.8}, $\stab_{\aut(V_{E_8})}(T_{16})/\pstab_{\aut(V_{E_8})}(T_{16}) \cong \aut(C)\cong\Sym_{8}\wr 2$ and thus $\stab_{\aut(V_{E_8})}(T_{16})$ has the shape $2^{2+12}.(\Sym_{8}\wr 2)$.

\medskip
\noindent \textbf{Case 3. ($\calC\cong 4^3\cdot 2^2$)}
In this case, one can see that $\calC_0= \Span_{\Z_2}\{ d(\Z_2^4), (10101010)\}$, $\calC_1= \Span_{\Z_2}\{ (1111\,0000), (0000\,1111), (1100\, 1100)\}$ and $\aut(\calC_0)\cong 2\wr \Sym_4$. The structure codes $(C,D)$ for $V_{A_4(\calC)}\cong V_{E_8}$ are given by
\[
C= \Span_{\Z_2}\{ \mathcal{E}_4^4, (10^3 10^3 10^3 10^3)\},\qquad D= \Span_{\Z_2}\{ (1^8 0^8), (0^8 1^8), (1^40^41^40^4) \}.
\]
Then, $\aut(C)\cong \aut(D)= \Sym_4\wr \Sym_4$ and $$P=\{\al\in C|\, \al\cdot\be\in C \text{ for all } \be\in D\}= \mathcal{E}_4^4,$$ which has dimension $12$. Thus,  $\pstab_{\aut(V_{E_8})}(T_{16})$ has the shape $2^{3+9}$.

Since $\aut(\calC)\cong 2^4.(2\wr\Sym_4)$, we have $\overline{\aut(\calC)}=2\wr \Sym_4$. Hence, we have $\overline{\aut(\calC)}\cong \aut(\calC_0)$ and $\stab_{\aut(V_{E_8})}(T_{16})/\pstab_{\aut(V_{E_8})}(T_{16}) \cong \aut(C) \cong \Sym_{4}\wr \Sym_4$. By Corollary \ref{Prop7.8} $\stab_{\aut(V_{E_8})}(T_{16})$ has the shape $2^{3+9}.(\Sym_{4}\wr \Sym_4)$.

\medskip
\noindent \textbf{Case 4. ($\calC\cong 4^4$)}
In this case, $\calC_0\cong \calC_1\cong H_8$, the Hamming $[8,4,4]$ code. The structure codes for $V_{A_4(\calC)}$ are given by
\[
C=\Span_{\Z_2}\{ d(\Z_2^8), e(H_8)\}, \qquad D=d(H_8).
\]

It is easy to see that $\aut(\calC)\cong \aut(H_8)\cong {\rm AGL}(3,2)$.
In this case, $P=\{\al\in C|\, \al\cdot\be\in C \text{ for all } \be\in D\}= \Span_{\Z_2}\{ d(\Z_2^8), e(1^8)\} $, which has dimension $9$ and $\overline{\aut(\calC)}\cong {\rm AGL}(3,2)$. Since the minimum weight of $H_8$ is $4$, it is clear that
$|\mathcal{H}|=1$. By Theorem \ref{Thm7.6} $\stab_{\aut(V_{E_8})}$ has the shape $2^{4+5}. (2\wr {\rm AGL}(3,2))$

\begin{remark}
The frame stabilizers for all Virasoro frames of the VOA $V_{E_8}$ were computed in \cite{GH} using a different method.
\end{remark}

\section{Frame stabilizer of the frame of $\tilde{V}_L$ associated to a $4$-frame}
Let $L$ be an even unimodular lattice of rank $n$ having $4$-frame $F$.
Let $\tilde{V}_L=V_L^+\oplus V_L^{T,+}$ be the VOA obtained by a $\Z_2$-orbifold construction from $V_L$.
Then $\tilde{V}_L$ has the Virasoro frame $T_{2n}$ associated to $F$.
In this section, we study the subgroup ${K}=\stab_{\aut(\tilde{V}_L)}(T_{2n})/\pstab_{\aut(\tilde{V}_L)}(T_{2n})$ of $\Aut(C)$ when the minimum weight of $\calC_0$ is greater than or equal to $4$, where $\calC=L/F$ is a $\Z_4$-code.
Note that the structure codes $(C,D)$ of $\tilde{V}_L$ associated to $T_{2n}$ were described in Proposition \ref{Prop DGH}.

\subsection{Stabilizer of $d(\mathcal{E}_n)$}
In this subsection, we determine the stabilizer in ${K}$ of the subcode $d(\mathcal{E}_n)$ of $C$.

\begin{proposition}\label{L3}
Let $L$ be an even unimodular lattice of rank $n$ having $4$-frame $F$ and $\calC=L/F$ a $\Z_4$-code of length $n$.
Let $T_{2n}$ be the Virasoro frame of $\tilde{V}_L$ associated to $F$.
Then the stabilizer of $d(\mathcal{E}_n)$ in ${K}$ has the shape $2^{\dim \calC_0^\perp}.\overline{\aut(\calC)}$.
\end{proposition}

\begin{proof}
Let $S$ be the stabilizer of $d(\mathcal{E}_n)$ in ${K}$.

For $k=1, \dots,n$, define $t_k= \exp(\frac{\pi i \al_k(0)}{4})$.
Then $t_k e^\beta=(-1)^{(\alpha_k/4,\beta)}e^\beta$.
Hence $\Pi_{k=1}^nt_k^{s_k}\in\aut(V_L^+)$ if and only if $(\sum_{k=1}^ns_k\al_k/4,\beta)\in\Z$ for all $\beta\in L$, that is, $\sum_{k=1}^ns_k\al_k/4\in L^*=L$.
Hence $\Pi_{k=1}^nt_k^{s_k}\in\aut(V_L^+)$ if and only if $\varphi_4(s_1,s_2,\dots,s_n)\in\calC$.
Assume that $h=\Pi_{k=1}^nt_k^{s_k}\in\aut(V_L^+)$.
Since $h$ is also an automorphism of $V_L$, it preserves $V_L^-$.
Since $L$ is unimodular and $h\circ V_L^-\cong V_L^-$, we have $h\circ V_L^{T,+}\cong V_L^{T,+}$.
Note that $\tilde{V}_L$ is a simple current extension of $V_L^+$.
Hence $h$ lifts to an automorphism of $\tilde{V}_L$ (\cite[Theorem 3.3]{Sh2}).
Clearly $t_k$ maps $\om^{\pm}(\al_k) $ to $\om^{\mp}(\al_k)$, and fixes $\om^\pm(\al_i)$ for all $i\neq k$.
Since $t_k^2=1$ on $T_{2n}$, $h=\Pi_{k=1}^nt_k^{\varphi_2(s_k)}$ on $T_{2n}$.
Thus, $A=\{\Pi_{k=1}^n \bar{t}_k^{s_k}\mid \varphi_2(s_1,s_2,\dots,s_n)\in\calC_1=\calC_0^\perp\}$ is a subgroup of $S$ isomorphic to $2^{\dim\calC_0^\perp}$, where $\bar{}\ $ is the canonical homomorphism from $\stab_{\aut(\tilde{V}_L)}(T_{2n})$ to $K$.

Let $g\in\aut(\calC)$. Then by (\ref{Eq:A4}) $g$ induces an automorphism of $\aut(L)$ preserving $F$.
Then it lifts to an automorphism $\tilde{g}$ of $\aut(\tilde{V}_L)$ preserving $T_{2n}$ (cf.\ \cite[Corollary 10.4.8]{FLM}).
Since $R\cap\aut(\calC)$ is trivial on $T_{2n}$, $\tilde{g}$ acts on $T_{2n}$ as $g\cdot (R\cap\aut(\calC))\in\overline{\aut(\calC)}= \aut(\calC)/(R\cap\aut(\calC))$ on $\{\{\omega^\pm(\alpha_i)\}\mid i\in\Omega_n\}$.
Thus we obtain a subgroup $A.\overline{\aut(\calC)}$ of $S$.

Let $g\in \stab_{\aut(\tilde{V}_L)}(T_{2n})$ such that $g\cdot \pstab_{\aut(\tilde{V}_L)}(T_{2n})\in S$.
Then it preserves the subVOA $V_{d(\mathcal{E}_n)}$ which is the code VOA associated to $d(\mathcal{E}_n)$.
It is easy to see that $V_{d(\mathcal{E}_n)}= V_{F}^+$.
We view $\tilde{V}_L$ as a module for $V_{d(\mathcal{E}_n)}$.
Since $F$ is a decomposable lattice and $n>1$, it is not obtained by Construction B.
Hence the restriction of $g$ on $V_{d(\mathcal{E}_n)}$ is a lift of $\bar{g}\in\aut(F)$ (\cite[Proposition 3.16]{Sh2}).
This shows that $g$ preserves $V_L^+$, and hence $g$ is a lift of $\bar{g}\in\aut(L)$ preserving $F$.
Thus $\bar{g}$ induces an automorphism of $L/F\cong\calC$, modulo the action of $\pstab_{\aut(\tilde{V}_L)}(T_{2n}).A$, and hence $\bar{g}\in \aut(\calC)$.
Therefore $S\cong A.\overline{\aut(\calC)}$.
\end{proof}

\begin{corollary}\label{Cor8} Assume that the minimum weight of $\calC_0$ is greater than $4$.
Then the group ${K}=\stab_{\aut(\tilde{V}_L)}(T_{2n})/\pstab_{\aut(\tilde{V}_L)}(T_{2n})$ has the shape $2^{\dim \calC_0^\perp}.\overline{\aut(\calC)}$
\end{corollary}
\begin{proof} By the assumption, any weight $4$ codeword in $C$ belongs to $d(\mathcal{E}_n)$.
Clearly, $d(\mathcal{E}_n)$ is generated by weight $4$ codewords.
Hence any automorphism of $C$ preserves $d(\mathcal{E}_n)$, and so does ${K}$.
This corollary then follows from Proposition \ref{L3}.
\end{proof}

\subsection{Pseudo Golay codes and the moonshine VOA}
In this subsection, we shall study  certain Virasoro frames of the moonshine VOA
$V^\natural$ arisen from pseudo Golay codes. The corresponding frame
stabilizers will also be computed.

Let $\calC$ be an extremal Type II $\Z_4$-code of length $24$ such that $\varphi_2(\calC)\cong G_{24}$.
Then the minimum Euclidean weight of $\calC$ is $16$, and
\[
A_4(\calC) =\frac{1}2 \{ (x_1, \dots, x_{24})\in \Z^{24}| \varphi_4(x_1, \dots,x_{24})\in \calC\}
\]
is isomorphic to the Leech lattice $\Lambda$.
By Proposition \ref{Prop DGH}, the structure codes for $\tilde{V}_{A_4(\calC)}=V_{A_4(\calC)}^+ \oplus V_{A_4(\calC)}^{T,+}\cong V^\natural$ are
\[
 C=\Span_{\Z_2}\{ d(\mathcal{E}_{24}), e(G_{24})\},  \qquad  D=\Span_{\Z_2}\{ d(G_{24}), e((1^{24}))\}.
\]
In this case, $
 P=\{\al\in C|\, \al\cdot\be\in C \text{ for all } \be\in D\}$ is equal to $D$ and has dimension $13$. Thus $\pstab_{\aut V^\natural}(T_{48}) =2^{13}$.
Since the minimum weight of $\calC_0\cong G_{24}$ is $8$, we obtain the following by Corollary \ref{Cor8}.

\begin{theorem}
Let $\calC$  be an extremal Type II $\Z_4$-code of length $24$ such that $\varphi_2(\calC)\cong G_{24}$ and let $T_{48}$ be the Virasoro frame of
$V^\natural= V_{A_4(\calC)}^+ \oplus V_{A_4(\calC)}^{T,+}$ associated to $\calC$.
Then,
\[
\stab_{\aut(V^\natural)}(T_{48})/\pstab_{\aut(V^\natural)}(T_{48})\cong 2^{\dim\calC_0^\perp}.\overline{\aut(\calC)}.\
\]
\end{theorem}

\begin{remark} In Rains \cite{R}, $13$ non-isomorphic extremal Type II $\Z_4$-codes satisfying $\varphi_2(\calC)\cong G_{24}$ were given as pseudo Golay codes. It was checked by Masaaki Harada and Akihiro Munemasa \cite{HM} that any extremal Type II $\Z_4$-codes satisfying $\varphi_2(\calC)\cong G_{24}$ is a pseudo Golay code.
\end{remark}

Next, we shall study  few examples from \cite[Fig. 2]{R}.

\begin{example}
 $\calC$ is generated by
\[ \small
\begin{pmatrix}
1   0   0   0   0   0   0   0   0   0   0   0   1   3   0   0   2   1   1   1
0   1   2   3\\
0   1   0   0   0   0   0   0   0   0   0   0   1   2   1   0   2   3   0   0
3   1   1   1\\
0   0   1   0   0   0   0   0   0   0   0   0   3   3   1   1   2   3   3   0
1   2   0   0\\
0   0   0   1   0   0   0   0   0   0   0   0   0   3   3   1   1   2   3   3
0   1   2   0\\
0   0   0   0   1   0   0   0   0   0   0   0   0   0   3   3   1   1   2   3
3   0   1   2\\
0   0   0   0   0   1   0   0   0   0   0   0   2   0   2   3   3   3   1   2
1   1   2   3\\
0   0   0   0   0   0   1   0   0   0   0   0   1   0   1   2   1   2   3   1
1   2   2   3\\
0   0   0   0   0   0   0   1   0   0   0   0   1   3   1   1   0   0   2   3
0   2   3   3\\
0   0   0   0   0   0   0   0   1   0   0   0   1   3   0   1   3   3   0   2
2   1   3   0\\
0   0   0   0   0   0   0   0   0   1   0   0   0   1   3   0   1   3   3   0
2   2   1   3\\
0   0   0   0   0   0   0   0   0   0   1   0   1   2   2   3   2   0   3   3
3   3   3   2\\
0   0   0   0   0   0   0   0   0   0   0   1   2   1   0   2   3   0   0   3
1   1   1   1
\end{pmatrix}_.
\]
In this case, $|\aut(\calC)|=12144$, $\aut(\calC)\cong {\rm SL}_2(23)$ and $\overline{\aut(\calC)}\cong {\rm PSL}_2(23)$.
Hence $\stab_{\aut(V^\natural)}(T_{48})$ has the shape $2^{13}. (2^{12}.{\rm PSL}_2(23))$.
Recall that the subgroup $H=\{\exp(\pi iv(0))\mid v\in \Lambda\}$ of $\Aut(V^\natural)$ is isomorphic to an extraspecial $2$-group $2^{1+24}_+$ (cf. \cite[(10.4.51)]{FLM}).
Since for any Ising vector $e$ in $T_{48}$, $\tau_{e}\in H$, $\pstab_{\Aut(V^\natural)}(T_{48})(\cong 2^{13})$ is a subgroup of $H$.
By Proposition \ref{L3}, we obtain a subgroup $\pstab_{\Aut(V^\natural)}(T_{48}).A\cong 2^{13}.2^{12}$ of $H$.
Comparing the orders, we obtain $\pstab_{\Aut(V^\natural)}(T_{48}).A=H\cong2^{1+24}_+$, and $\stab_{\aut(V^\natural)}(T_{48})$ has the shape $2^{1+24}_+.{\rm PSL}_2(23)$.
\end{example}

\begin{example}
 $\calC$ is generated by
 \[ \small
\begin{pmatrix}
1   0   0   0   0   0   0   0   0   0   0   0   222201113131\\
0   1   0   0   0   0   0   0   0   0   0   0   203131231223\\
0   0   1   0   0   0   0   0   0   0   0   0   301212231112\\
0   0   0   1   0   0   0   0   0   0   0   0   010330223311\\
0   0   0   0   1   0   0   0   0   0   0   0   322132330121\\
0   0   0   0   0   1   0   0   0   0   0   0   313332121022\\
0   0   0   0   0   0   1   0   0   0   0   0   132033122312\\
0   0   0   0   0   0   0   1   0   0   0   0   211012130031\\
0   0   0   0   0   0   0   0   1   0   0   0   332301233322\\
0   0   0   0   0   0   0   0   0   1   0   0   303323312030\\
0   0   0   0   0   0   0   0   0   0   1   0   133021001211\\
0   0   0   0   0   0   0   0   0   0   0   1   213303320101
\end{pmatrix}_.
\]
Then $|\aut(\calC)|=6$, $\aut(\calC)\cong \Z_2\times\Z_3$ and $\overline{\aut(\calC)}\cong \Z_3$.
Thus, $\stab_{\aut(V^\natural)}(T_{48})$ has the shape $2^{13}. (2^{12}.3)$.
By the same argument in the previous example, the shape of $\stab_{\aut(V^\natural)}(T_{48})$ is $2^{1+24}_+.3$.
\end{example}

\begin{remark}
 By the examples above, we note that $\stab_{\aut(V)}(T_{48})/{\pstab}_{\aut(V)}(T_{48})$ may be strictly smaller than $\aut(C)$ in general.  In fact, $\stab_{\aut(V^\natural)}(T_{48})/{\pstab}_{\aut(V^\natural)}(T_{48})$ can be quite small compare to $\aut(C)\,(\cong 2\wr M_{24})$.
\end{remark}

\subsection{Minimum weight of $\calC_0$ is $4$}
Next we shall consider the case where the minimum weight of $\calC_0$ is $4$.
Note that the case where the minimum weight of $\calC_0$ is greater than $4$ was done in Corollary \ref{Cor8}.
Remark that $n\in8\Z$ since $A_4(\mathcal{C})$ is even unimodular.

First, we recall the following easy lemmas.

\begin{lemma}\label{LemEn}
Let $W$ be a binary code of length $2n$ isomorphic to $d(\mathcal{E}_{n})$.
Then the following hold.
\begin{enumerate}
\item $W$ is generated by $W(4)$.
\item $|w_1\cap w_2|\in2\Z$  for all $w_1, w_2\in W$.
\item For any $w\in W(4)$, $|\{w^\prime\in W(4)\mid |w\cap w^\prime|=2\}|=2n-4$.
\item Let $k\ge2$ and let $w_1,w_2,\dots,w_k\in W(4)$ such that $|w_i\cap w_j|=2\delta_{|i-j|,1}$ for $i\neq j$.
Then $|\{w\in W(4)\mid |w\cap w_i|=2\delta_{i,k}, 1\le i\le k\}|=n-k-1$.
\item If $n>8$ then $W^\perp(4)=W(4)$.
\item If $n>2$ then for any $w\in W(4)$ there exist unique $w_1,w_2\in W^\perp(2)$ such that $w=w_1+w_2$.
\item If $n>2$ then for any distinct $w_1,w_2\in W^\perp(2)$, $w_1+w_2\in W(4)$.
\end{enumerate}
\end{lemma}


\begin{lemma}\label{LemC} Assume that the minimum weight of $\calC_0$ is $4$.
Then the following hold.
\begin{enumerate}
\item $C(4)=\{d(x)\mid x\in \mathcal{E}_n(2)\}\cup
\{ e(y)+d(z)\mid y\in\calC_0(4),\ z\in\mathcal{E}_n,\ z\subset y\}.$
\item For $e(y)+d(z),d(x)\in C(4)$, $|(e(y)+d(z))\cap d(x)|=0$ if and only if $|x\cap y|=0$.
Moreover $|(e(y)+d(z))\cap d(x)|=2$ if and only if $x\subset y$.
\item Let $u=e(y)+d(z),u^\prime=e(y^\prime)+d(z^\prime)\in C(4)$, where $y,y^\prime\in\calC_0(4)$.
If $|u\cap u^\prime|=2$ then $|y\cap y^\prime|\in\{2,4\}$.
\end{enumerate}
\end{lemma}
\begin{proof} The lemma  can be deduced easily from the fact that $C=\Span_{\Z_2}\{ d(\mathcal{E}_n),e(\calC_0)\}$ (Proposition \ref{Prop DGH}).
\end{proof}

\begin{lemma}\label{L7-1} Let $s_1,s_2,\dots, s_{n/2-1}$ be weight $4$ elements in $\Z_2^n$ such that $|s_i\cap s_j|=2\delta_{|i-j|,1}$ if $i\neq j$.
Let $t_1,t_2,\dots,t_{n/2}$ be weight $2$ elements in $\Z_2^{n}$ such that $s_i=t_i+t_{i+1}$.
Then $\Span_{\Z_2}\{ d(t_i), e(s_j)\mid 1\le i\le n/2, 1\le j\le n/2-1\}$ is isomorphic to $d(\mathcal{E}_n)$.
\end{lemma}
\begin{proof} Set $r_{2i-1}=d(t_i)$, $r_{2j}=e(s_{j})+d(t_{j+1})$ for $1\le i\le n/2$, $1\le j\le n/2-1$.
Then the weight of $r_i$ is $4$, and $|r_i\cap r_j|=2\delta_{|i-j|,1}$ if $i\neq j$.
Hence, $\Span_{\Z_2}\{ d(t_i), e(s_j)\mid 1\le i\le n/2, 1\le j\le n/2-1\}=\Span_{\Z_2}\{ r_1,r_2,\dots, r_{n-1}\}\cong d(\mathcal{E}_n)$.
\end{proof}

Let $E$ be a subcode of $C$ isomorphic to $d(\mathcal{E}_n)$ such that $E\neq d(\mathcal{E}_n)$.
Set $Y=\{ c\in\calC_0\mid (e(c)+d(\mathcal{E}_n))\cap E\neq\emptyset\}$.
In order to determine $E$, we need some lemmas.

\begin{lemma}\label{L7-00} There exist $y_1,y_2,\dots, y_{n/2-1}\in Y(4)$ and $u_1,u_2,\dots, u_{n/2-1}\in E(4)$ such that $u_i\in (e(y_i)+d(\mathcal{E}_n))\cap E(4)$ and $|y_i\cap y_j|=|u_i\cap u_j|=2\delta_{|i-j|,1}$ if $i\neq j$.
\end{lemma}
\begin{proof} By Lemma \ref{LemEn} (1) and the assumption $E\neq d(\mathcal{E}_n)$, there is $u_1\in E(4)\setminus d(\mathcal{E}_n)$.
By Lemma \ref{LemC} (1), $u_1=e(y_1)+d(z_1)$ for some $y_1\in Y(4)$ and $z_1\in\mathcal{E}_n$ satisfying $z_1\subset y_1$.
By Lemma \ref{LemEn} (3), $|\{u\in E(4)\mid |u\cap u_1|=2\}|=2n-4$.
By Lemma \ref{LemC} (2), $$|\{d(x)\in E(4)\mid x\in\mathcal{E}_n(2),\ |d(x)\cap u_1|=2\}|\le \binom{4}{2}=6.$$
Clearly $$\left|\{u=e(y_1)+d(z)\in E(4)\mid z\in\mathcal{E}_n,\ z\subset y_1,\ |u\cap u_1|=2\}\right|\le \binom{4}{2}=6.$$
If $n=8$ and if the both equalities hold then $E$ contains a subcode isomorphic to the extended Hamming code of length $8$, which contradicts $E\cong d(\mathcal{E}_8)$.
If $n>8$, then $2n-4> 6+6=12$.
Hence there exists $u_2\in E(4)$ such that $|u_1\cap u_2|=2$ and $u_2=e(y_2)+d(z_2)$ for some $y_2\in Y(4)$ satisfying $y_2\neq y_1$.
Then by Lemma \ref{LemC} (3), $|y_1\cap y_2|=2$.

Let $k\in\Z$ satisfying $2\le k\le n/2-2$.
Assume that there exist $u_1=e(y_1)+d(z_1),u_2=e(y_2)+d(z_2),\dots, u_k=e(y_k)+d(z_k)\in E(4)$, $y_i\in Y$ such that $|u_i\cap u_j|=|y_i\cap y_j|=2\delta_{|i-j|,1}$ if $i\neq j$.
Let us show that there exist $u_{k+1}=e(y_{k+1})+d(z_{k+1})\in E(4)$ and $y_{k+1}\in Y(4)$ such that $|y_{k+1}\cap y_j|=|u_{k+1}\cap u_j|=2\delta_{j,k}$ for $1\le \forall j\le k$.

By Lemma \ref{LemEn} (4), $|\{u\in E(4)\mid |u\cap u_j|=2\delta_{j,k}\}|=n-k-1$.
By $|u_i\cap u_j|=2\delta_{|i-j|,1}$ and Lemma \ref{LemC} (2), $$|\{d(x)\in (d(\mathcal{E}_n)\cap E)(4)\mid |d(x)\cap u_j|=2\delta_{j,k}\}|\le 1.$$
Set $x_1=y_1\setminus y_2$, $x_i=y_i\cap y_{i-1}$ ($2\le i\le k$) and set $I= \Span_{\Z_2}\{ u_i,d(\mathcal{E}_n)\mid 1\le i\le k\}$.
Let $u\in I(4)\setminus d(\mathcal{E}_n)$ such that $|u\cap u_j|=2\delta_{j,k}$.
Then $u=\sum_{i=p}^k u_i+d(x_p)$ for some $1\le p\le k$, and we obtain
$$|\{u\in I(4)\setminus d(\mathcal{E}_n)\mid |u\cap u_j|=2\delta_{j,k}\}|\le k.$$
It follows from $k\le n/2-2$ that $n-k-1\ge k+3>k+1$.
Hence, there exists $u_{k+1}\in E(4)\setminus I$ such that $|u_{k+1}\cap u_j|=2\delta_{j,k}$.
By Lemma \ref{LemC} (1), there exists $y_{k+1}\in Y(4)$ such that $u_{k+1}=e(y_{k+1})+d(z_{k+1})$ for some $z_{k+1}\in\mathcal{E}_n$ and $z_{k+1}\subset y_{k+1}$.
Since $u_{k+1}\notin I$, we have $y_{k+1}\neq y_k$.
Hence $|y_{k+1}\cap y_k|=2$ by Lemma \ref{LemC} (3).
Let $q\in\Z$ such that $1\le q\le k-1$.
Clearly $|y_{k+1}\cap y_q|\neq 3,4$.
If $|y_{k+1}\cap y_q|=1$ then the weight of $\sum_{i=q}^{k+1}y_i$ is $2$, which contradicts that the minimum weight of $\calC_0$ is $4$.
Since $y_{k+1}\notin I$, we have $|y_{k+1}\cap y_q|\neq2$.
Hence $|y_{k+1}\cap y_q|=0$, and we obtain desired elements $y_{k+1}\in Y(4)$ and $u_{k+1}\in E(4)$.
Thus by induction, we obtain this lemma.
\end{proof}

Let $x_1,x_2,\dots, x_{n/2}\in \mathcal{E}_n(2)$ such that $y_i=x_i+x_{i+1}$.
Then $x_i=y_i\cap y_{i-1}$ ($2\le i\le n/2-1$), $\sum_{i=1}^{n/2}x_i=\Omega_n$ and $|x_i\cap y_j|=2\delta_{i,j}+2\delta_{i,j-1}$.
Set $X=\Span_{\Z_2}\{ x_1,x_2,\dots,x_{n/2}\}$, $\tilde{Y}=\Span_{\Z_2}\{ y_1,y_2,\dots,y_{n/2-1}\}$ and $U=\Span_{\Z_2}\{ u_1,u_2,\dots, u_{n/2-1}\}$.

\begin{lemma}\label{L7-01} Let $u=e(y)+d(z)\in E(4)$.
Then $y\in \tilde{Y}(4)$.
\end{lemma}
\begin{proof} Let $y^\prime\in \tilde{Y}(4)$.
Then there exists $u^\prime=e(y^\prime)+d(z^\prime)\in E(4)$.
Since the minimum weight of $\calC_0$ is $4$, $|y\cap y^\prime|\in\{0,1,2,4\}$.

First, we will show that $|y\cap y^\prime|\in2\Z$.
Suppose that $|y\cap y^\prime|=1$.
Then $|y\cap x|\in\{0,1\}$ for all $x\in X(2)$ since the minimum weight of $\calC_0$ is $4$.
Let $\{x^1,x^2,x^3,x^4\}\subset X(2)$ such that $|y\cap x_i|=1$.
Now, we view $\Span_{\Z_2}\{ x^i+x^j,y\mid 1\le i,j\le 4\}$ as a code of length $8$.
Then it is isomorphic to the extended Hamming code of length $8$.
Up to coordinates, we may assume that
\begin{eqnarray*}
x^1&=&(11000000),\\
x^2&=&(00110000),\\
x^3&=&(00001100),\\
x^4&=&(00000011),\\
y&=&(01010101).
\end{eqnarray*}
If $n=8$ then $y^\prime\in\Span_{\Z_2}\{x^i+x^j\mid 1\le i,j\le 4\}$, and hence $|y\cap y^\prime|\in2\Z$, which contradicts $|y\cap y^\prime|=1$.
Hence $n>8$.
Let $u^1,u^2,u^3\in U(4)$ such that $u^i\in e(x^i+x^{i+1})+d(\mathcal{E}_n)$.
Then $|u^i\cap u^j|=2\delta_{i+1,j}$ if $i\neq j$.
We may also view $u^i$, $i=1,2,3$, as codewords of length $16$.
Then, up to coordinates, we may assume that
\begin{eqnarray*}
u^1&=&(0101010100000000),\\
u^2&=&(0000010101010000),\\
u^3&=&(0000000001010101).
\end{eqnarray*}
If $|u\cap u^i|=2$ then there exists $w\in U$ such that $|u\cap w|=1$ since $n>8$.
This contradicts Lemma \ref{LemEn} (2).
Hence $|u\cap u^i|=0$.
Since $u\notin \Span_{\Z_2}\{u^1, u^2,u^3\}$, up to coordinates, we may assume that \begin{eqnarray*}u&=&(0010001000100010). \end{eqnarray*}
By Lemma \ref{LemEn} (6), $E^\perp$ contains
\begin{eqnarray*}
s^1&=&(0101000000000000),\\
s^2&=&(0000010100000000),\\
s^3&=&(0000000001010000),\\
s^4&=&(0000000000000101).
\end{eqnarray*}
Moreover, $E^\perp$ contains
$$s^5=(0010001000000000),\quad (0010000000100000),\quad {\rm or}\quad (0010000000000010).$$
If $s^5=(0010001000000000)\in E^\perp$ then by Lemma \ref{LemEn} (7),
$$s^1+s^5=(0111000100000000)\in E(4),$$ which contradicts the minimum weight of $\calC_0$ is $4$.
We may obtain contradictions by similar arguments for $s^5=(0010000000100000)\text{ and  } (0010000000000010).$
Thus $|y\cap y^\prime|\neq1$.

Therefore, we have $|y\cap y^\prime|\in2\Z$, and $y\in\tilde{Y}^\perp(2)$.
If $n>8$ then by Lemma \ref{LemEn} (5), we obtain $y\in \tilde{Y}(4)$.
If $n=8$ and $y\notin \tilde{Y}$ then $\Span_{\Z_2}\{ \tilde{Y},y\}$ is isomorphic to the extended Hamming code of length $8$.
By the arguments above, up to coordinates, $E$ contains
\begin{eqnarray*}
u^1&=&(0101010100000000),\\
u^2&=&(0000010101010000),\\
u^3&=&(0000000001010101),\\
u&=&(0001000100010001).
\end{eqnarray*}
However, these generate the extended Hamming code of length $8$, which contradicts $E\cong d(\mathcal{E}_n)$.
Therefore $y\in \tilde{Y}(4)$.
\end{proof}

\begin{proposition}\label{L7-0} Let $E$ be a subcode of $C$ isomorphic to $d(\mathcal{E}_n)$ such that $E\neq d(\mathcal{E}_n)$.
Set $Y=\{ c\in\calC_0\mid (e(c)+d(\mathcal{E}_n))\cap E\neq\emptyset\}$.
Let $y_i\in Y(4)$ and $u_i\in E(4)$ be codewords given in Lemma \ref{L7-00}.
Let $x_1,x_2,\dots,x_{n/2}\in \mathcal{E}_n(2)$ such that $y_i=x_i+x_{i+1}\in Y$.
Let $w_1,w_2,\dots, w_{n/2}\in\mathcal{E}_n(2)$ satisfying $|x_i\cap w_j|=\delta_{i,j}+\delta_{i+1,j}$.
Then one of the following holds:
\begin{itemize}
\item $E=\Span_{\Z_2}\{ d(x_i), e(y_j)\mid 1\le i\le n/2, 1\le j\le n/2-1\}$.
\item $E=\Span_{\Z_2}\{ d(x_i), e(y_j)+d(w_j)\mid 1\le i\le n/2, 1\le j\le n/2-1\}$.
\end{itemize}
\end{proposition}
\begin{proof}
Set $\tilde{Y}=\Span_{\Z_2}\{ y_1,y_2,\dots, y_{n/2-1}\}$.
Let $u=e(y)+d(z)\in E(4)$, where $y\in {Y}(4)$ and $z\in d(\mathcal{E}_n)$.
Then by Lemma \ref{L7-01}, $y\in \tilde{Y}(4)$.

Set $U=\Span_{\Z_2}\{ u_1,u_2,\dots, u_{n/2 -1}\}$ and $X=\Span_{\Z_2}\{x_1,x_2,\dots,x_{n/2} \}$.
Let $d(x)\in E(4)$.
Then $\langle e(y^\prime)+d(z^\prime),d(x)\rangle=\langle e(y), d(x)\rangle=\langle y,x\rangle=0$ for any $e(y^\prime)+d(z^\prime)\in E$.
Hence $x\in \tilde{Y}^\perp(2)$.
Since $n\ge 2$, $\tilde{Y}^\perp(2)=X(2)$.
Thus by Lemma \ref{LemC} (1), we obtain
\begin{eqnarray}
\hspace{1cm}E(4)\subset\{d(x)\mid x\in X(2)\}\cup\{ u+d(x)\mid u=e(y)+d(z)\in U(4),x\in X, x\subset y\}.\label{E4}
\end{eqnarray}
It is easy to see that $|E(4)|=n\times (n-1)/2$ and the cardinality of the right hand in (\ref{E4}) is equal to $n/2+(n/2)(n/2-1)/2\times 4=n\times (n-1)/2$.
Hence the equality holds in (\ref{E4}) and $Y=\tilde{Y}$.
Moreover $u_i\in e(y_i)+\Span_{\Z_2}\{ d(x_i),d(x_{i+1})\}$ or $u_i\in e(y_i)+d(w_i)+\Span_{\Z_2}\{ d(x_i),d(x_{i+1})\}$, where $w_i\in\mathcal{E}_n(2)$ such that $|w_i\cap x_i|=|w_i\cap x_{i+1}|=1$.

Let us determine $E$.
If $u_i\in e(y_i)+\Span_{\Z_2}\{ d(x_i),d(x_{i+1})\}$ then $e(y_i)\in E$.
Suppose $u_{i+1}\in e(y_{i+1})+d(w_{i+1})+\Span_{\Z_2}\{d(x_{i+1}),d(x_{i+2})\}$.
Then $\langle e(y_i),u_{i+1}\rangle=1$, which contradicts Lemma \ref{LemEn} (2).
Hence $u_{i+1}\in e(y_{i+1})+\Span_{\Z_2}\{d(x_{i+1}),d(x_{i+2})\}$, and $e(y_{i+1})\in E$.
Thus $E=\Span_{\Z_2}\{ d(x_i), e(y_j)\mid 1\le i\le n/2, 1\le j\le n/2-1\}$.

By the similar arguments, if $u_i\in e(y_i)+d(w_i)+\Span_{\Z_2}\{ d(x_i),d(x_{i+1})\}$, then $e(y_j)+d(w_j)\in E$ for all $j$.
Hence $E=\Span_{\Z_2}\{ d(x_i), e(y_j)+d(w_j)\mid 1\le i\le n/2, 1\le j\le n/2-1\}$.
\end{proof}

\begin{lemma}\label{L7-2} Let $L$ be an even unimodular lattice of rank $n$ having $4$-frame $F$.
Let $\calC=L/F$ be a type II $\Z_4$-code.
Assume that $\calC_0$ contains a subcode $Y$ isomorphic to $d(\mathcal{E}_{n/2})$ and that the minimum weight of $\calC_0$ is $4$.
Then there exists $\{f_i\mid 1\le i\le n\}\subset\R^n$ such that $(f_i,f_j)=2\delta_{ij}$, $f_i+f_j\in L$, $f_i\notin L$, and $\{f_{2k-1}\pm f_{2k}\mid 1\le k\le n/2\}=F$.
\end{lemma}
\begin{proof} There exist $x_1,x_2,\dots,x_{n/2}\in\mathcal{E}_n(2)$ such that $x_i+x_j\in Y(4)$.
By the definition of $\calC$, $$L\cong A_4(\calC)=\frac{1}{2}\{(c_1,c_2,\dots,c_n)\in\Z^n\mid \varphi_4(c_1,c_2,\dots,c_n)\in\calC\}.$$
It is easy to see that $\frac{1}{2}\{c\in \Z^n\mid \varphi_2(c)=x_i\}\subset\R^n$ contains exactly $4$ vectors of norm $2$.
Hence, we can choose vectors $f_{2k-1}, f_{2k}$ of norm $2$ satisfying $f_{2k-1}\pm f_{2k}\in F$.
Then $f_i+f_j\in L$, $(f_{2k-1},f_{2k})=0$, and $\{f_i\mid 1\le i\le n\}$ is an orthogonal basis of $\R^n$.
By the assumption, there are no weight $2$ elements in $\calC_0$, and hence $f_i\notin L$.
\end{proof}

Now we recall a relation between codes and lattices.

\begin{proposition}\cite[Theorem 1]{KKM}\label{KKM} Let $L$ be an even unimodular lattice of rank $n$.
Assume that there exists an orthogonal basis $\{f_i\mid 1\le i\le n\}$ of norm $2$ of $\R^n$ such that $f_i+f_j\in L$ and $f_i\notin L$.
Then there exists a binary Type II code $W$ of length $n$ such that
$$L=\sum_{i,j\in\Omega_n}\Z(f_i+f_j)+\frac{1}{2}\sum_{w\in W}\Z f_w+\Z(\frac{1}{4}f_{\Omega_n}-\varepsilon(n) f_1),$$ where $f_w=\sum_{i\in w} f_i$ and $\varepsilon(n)=1$ if $n\in8+16\Z$ and $\varepsilon(n)=0$ if $n\in16\Z$.
\end{proposition}

\begin{proposition}\cite{FLM}\label{P7-1}
Let $L$ be an even unimodular lattice of rank $n$.
Let $\{f_i\mid 1\le i\le n\}$ be an orthogonal basis of $\R^n$ of norm $2$ satisfying $f_i+f_j\in L$ and $f_i\notin L$.
Then there exists an automorphism $\sigma$ of $\aut(\tilde{V}_L)$ such that $\sigma$ acts on the Virasoro frame associated to $4$-frame $\{f_{2i-1}\pm f_{2i}\mid 1\le i\le n/2\}$ of $L$ as
\begin{eqnarray}
\sigma (\omega^+(f_{2i-1}+f_{2i}))=\omega^+(f_{2i-1}+f_{2i}),\notag\\
 \sigma(\omega^-(f_{2i-1}+f_{2i}))=\omega^+(f_{2i-1}-f_{2i}),\label{Eq7-1}\\
 \sigma(\omega^+(f_{2i-1}-f_{2i}))=\omega^-(f_{2i-1}+f_{2i}),\notag\\
\sigma (\omega^-(f_{2i-1}-f_{2i}))=\omega^-(f_{2i-1}-f_{2i}).\notag
\end{eqnarray}
\end{proposition}
\begin{proof} By Proposition \ref{KKM} and \cite[(11.2.6)]{FLM}, we can obtain a triality automorphism $\sigma$ of $\tilde{V}_L$, which is similar to one of the moonshine VOA $V^\natural$.
The action of $\sigma$ on the Virasoro frame associated to $4$-frame $\{f_{2i-1}\pm f_{2i}\mid 1\le i\le n/2\}$ was described in \cite[Section 4]{DGH} (cf. \cite[Corollary 11.24]{FLM}).
\end{proof}

\begin{proposition}\label{P7-2} Let $L$ be an even unimodular lattice of rank $n$ and let $\calC=L/F$ be a type II $\Z_4$-code of length $n$.
Assume that the minimum weight of $\calC_0$ is $4$ and that $\calC_0$ contains a subcode $Y$ isomorphic to $d(\mathcal{E}_{n/2})$.
Let $\{y_1,y_2,\dots, y_{n/2-1}\}$ be a basis of $Y$ of weight $4$ such that $|y_i\cap y_j|=2\delta_{|i-j|,1}$ if $i\neq j$.
Let $x_1,x_2,\dots,x_{n/2}$ be codewords in $\mathcal{E}_n$ of weight $2$ such that $y_i=x_i+x_{i+1}$.
Then there exists an automorphism $g\in \stab_{\aut(\tilde{V}_L)}(T_{2n})/\pstab_{\aut(\tilde{V}_L)}(T_{2n})$ such that $g(d(\mathcal{E}_n))=\Span_{\Z_2}\{ e(y_i),d(x_j)\mid 1\le i\le n/2-1,\ 1\le j\le n/2\}$.
\end{proposition}
\begin{proof} Let $F=\{f_{2i-1}\pm f_{2i}\mid 1\le i\le n/2\}$ be the $4$-frame of $L$ obtained in the proof of Lemma \ref{L7-2}.
By Lemma \ref{L7-2} and Proposition \ref{P7-1}, there exists a triality automorphism $\sigma$ in $\stab_{\aut(\tilde{V}_L)}(T_{2n})$ satisfying (\ref{Eq7-1}).

Let $w_1,w_2,\dots,w_{n/2-1}\in \mathcal{E}_n(2)$ such that $|x_i\cap w_j|=\delta_{i,j}+\delta_{i+1,j}$.
Then $\{x_i,w_j\mid 1\le i\le n/2,\ 1\le j\le n/2-1\}$ is a basis of $\mathcal{E}_n$.
We view $\{(0^4),d(x_i)\}$ and $\{(0^4),d(w_i)\}$ as binary codes of length $4$.
Then $V_{\Z f_{2i-1}\oplus\Z f_{2i}}^+=V_{\{(0^4),d(x_i)\}}$, and $\sigma(d(x_i))=d(x_i)$.
However, $V_{\Z f_{2i}\oplus\Z f_{2i+1}}^+=V_{\{(0^4),d(w_i)\}}$ is not preserved by $\sigma$.
It follows from (\ref{Eq7-1}) that $\sigma(d(w_i))=e(y_i)+d(x_{i+1})$.
Hence $g=\bar{\sigma}\in \stab_{\aut(\tilde{V}_L)}(T_{2n})/\pstab_{\aut(\tilde{V}_L)}(T_{2n})$ sends $d(\mathcal{E}_n)$ to $\Span_{\Z_2}\{ e(y_i),d(x_j)\mid 1\le i\le n/2-1,\ 1\le j\le n/2\}$ (see the proof of Lemma \ref{L7-1}).
\end{proof}

\begin{lemma}\label{L7-4} Assume that $C$ contains $E=\Span_{\Z_2}\{ d(x_i), e(y_j)\mid 1\le i\le n/2, 1\le j\le n/2-1\}$ isomorphic to $d(\mathcal{E}_n)$.
Then there exists an automorphism $g\in\stab_{\Aut(\tilde{V}_L)}(T_{2n})$ such that $g(E)=\Span_{\Z_2}\{d(x_i), e(y_j)+d(w_j)\mid 1\le i\le n/2, 1\le j\le n/2-1\}\subset C$ (see Proposition \ref{L7-0} for the definition of $w_i$).
\end{lemma}
\begin{proof} By Lemma \ref{L7-2}, we obtain the orthogonal basis $\{f_i\}$ of norm $2$ of $\R^n$ such that $f_i+f_j\in L$, $f_i\notin L$.
By Proposition \ref{KKM}, we have $v=\frac{1}{4}f_{\Omega_n}-\varepsilon(n)f_1\in L$.
Then one can see that $\exp(\pi iv(0))\in\Aut(\tilde{V_L})$ is a desired automorphism.
\end{proof}

Let ${\mathcal{H}}$ be the set of all subcodes of $C$ isomorphic to $d(\mathcal{E}_n)$.
By Propositions \ref{L7-0} and \ref{P7-2} and Lemma \ref{L7-4}, we obtain the following theorem.

\begin{theorem}\label{MT} Let $L$ be an even unimodular lattice of rank $n$ having $4$-frame $F$.
Set $\calC=L/F$.
Assume that the minimum weight of $\calC_0$ is $4$.
Let $T_{2n}$ be the Virasoro frame of $\tilde{V}_L$ associated to $F$.
Then $\stab_{\aut(\tilde{V}_L)}(T_{2n})/\pstab_{\aut(\tilde{V}_L)}(T_{2n})$ is transitive on ${\mathcal{H}}$, and it is generated by the subgroup of shape $2^{\dim\calC_0^\perp}.\overline{\aut(\calC)}$ and the triality automorphisms in Proposition \ref{P7-1}.
Moreover, $|\stab_{\aut(\tilde{V}_L)}(T_{2n})/\pstab_{\aut(\tilde{V}_L)}(T_{2n}):2^{\dim\calC_0^\perp}.\overline{\aut(\calC)}|=|{\mathcal{H}}|$.
\end{theorem}

Next let us consider the stabilizer of $d(\mathcal{E}_n)$ in $\Aut(C)$.

\begin{lemma} Assume that the minimum weight of $\calC_0$ is $4$.
Then the stabilizer of $d(\mathcal{E}_n)$ in $\aut(C)$ has the shape $2^{\dim \calC_0^\perp}:\aut(\calC_0)$.
\end{lemma}
\begin{proof} Clearly, $2^{\dim \calC_0^\perp}:\aut(\calC_0)$ is a subgroup of $\aut(d(\mathcal{E}_n))$.
Let $g\in \aut(d(\mathcal{E}_n))$.
Then $g\in\aut(d(\mathcal{E}_n)^\perp)$.
It is easy to see that $d(\mathcal{E}_n)^\perp=\Span_{\Z_2}\{ d(\Z^n),e(1^n)\}$ and that $d(\mathcal{E}_n)^\perp(2)=\{d(x)\mid x\in\Z_2^n(1)\}$ since $n>2$.
Hence $g\in 2\wr\Sym_n$, and $g\in 2\wr\aut(\calC_0)$.
It is easy to see that the subgroup of $2^n$ preserving $C$ is isomorphic to $2^{\dim \calC_0^\perp}$.
\end{proof}

As a corollary, we have the following proposition.

\begin{corollary}\label{index:Vtilde}
Set ${K}= \stab_{\aut(\tilde{V}_L)}(T_{2n})/\pstab_{\aut(\tilde{V}_L)}(T_{2n})$ and assume that the minimum weight of $\calC_0$ is $4$.
Then
\[
 |\aut(C): {K}| = |\aut(\calC_0): \overline{\aut(\calC)}|.
\]

\end{corollary}

\begin{remark} Let $g\in \Aut(C)$.
Then there exists $h\in {K}$ such that $hg$ fixes $d(\mathcal{E}_n)$ by Theorem \ref{MT}.
Hence $g\in {K}$ if and only if $hg$ belongs to the stabilizer of $d(\mathcal{E}_n)$ in ${K}$ determined in Proposition \ref{L3}.
\end{remark}

\subsection{$\Z_2$-orbifold construction of $V_{E_8}$}
In this section, we consider the $\Z_2$-orbifold construction of $V_{E_8}$.

Let $\calC$ be the Type II $\Z_4$-code of length $8$ generated by
\[
\begin{pmatrix}
3111\, 3111\\
1111\, 2000\\
1320\, 1100\\
1010\, 1032
\end{pmatrix}_.
\]
Then the minimum weight of $\calC_0$ is $4$ and $A_4(\calC)\cong E_8$.
Then $\tilde{V}_{A_4(\calC)}= V_{A_4(\calC)}^+\oplus V_{A_4(\calC)}^{T,+}\cong V_{E_8}$.
Let $T_{16}$ be the Virasoro frame associated to $\calC_0$.
In this case,  $\calC_0\cong \calC_1\cong H_8$ and the structure codes for $\tilde{V}_{A_4(\calC)}\cong V_{E_8}$ are given by
\[
C=\Span_{\Z_2}\{ d(\mathcal{E}_8), e(H_8)\}\cong \RM(2,4), \qquad  D=\Span_{\Z_2}\{ d(H_8), e(1^8)\} \cong \RM(1,4).
\]
The code $P=\{\al\in C|\, \al\cdot\be\in C \text{ for all } \be\in D\}$ is equal to $D$ and has dimension $5$.

It is also well known that
\[
 \aut(C)\cong \aut(\mathrm{RM}(2,4)) \cong {\rm AGL}(4,2) \cong 2^4: {\rm GL}(4,2)
 \]
 and
\[
 \overline{\aut(\calC)}\cong \aut(\calC_0)\cong \aut(H_8)\cong {\rm AGL}(3,2)\cong 2^3:{\rm GL}(3,2).
\]

Since $\overline{\aut(\calC)}\cong \aut(\calC_0)$, by Corollary \ref{index:Vtilde}, we have $|\aut(C): {K}|=1$. Thus,
$$\stab_{\aut(\tilde{V}_{E_8})}(T_{16})/ \pstab_{\aut(\tilde{V}_{E_8})}(T_{16})\cong \aut(C)\cong 2^4:{\rm GL}(4,2)\cong {\rm AGL}(4,2)$$ and $\stab_{\aut(\tilde{V}_{E_8})}(T_{16})$ has the shape
$2^5.({\rm AGL}(4,2))$. \qed

\begin{remark}
The example above was also computed in \cite{GH}.
\end{remark}

\begin{remark}
Since the binary codes $\calC_0$ for the first three Type II $\Z_4$-codes $\calC$ in Example \ref{SsEx} contain weight $2$ codewords, we can not apply Corollary \ref{index:Vtilde} to the Virasoro frames of $\tilde{V}_{E_8}$ associated to $\calC$.
However, the structure code $C$ contains a subcode isomorphic to $d(\Z_2^8)$, and hence the Virasoro frame of $\tilde{V}_{E_8}$ associated to $\calC$ is conjugate to one of the frames of $V_{E_8}$ associated to $4$-frames of the lattice $E_8$ (cf.\ \cite[Section 5]{DGH}).
\end{remark}

\subsection{Frame stabilizer of the standard Virasoro frame of the moonshine
VOA}

Let us recall the standard construction of the Leech lattice from the binary Golay code $G_{24}$ of length $24$ (\cite[p131, Figure 4.12]{CS98}).
Viewing $0$ and $1$ as integers, the Leech lattice  $\Lambda$ is given by
$ \Lambda= \Lambda^0\cup \Lambda^1,$
where
\[
\begin{split}
\Lambda^0= &\frac{1}{\sqrt{8}} \{ 2c + 4x\ |\ c\in G_{24}, x\in \Z^{24}, \text{
and } \sum_{i=1}^{24} x_i\equiv 0\mod 2\},\\
\Lambda^1= &\frac{1}{\sqrt{8}} \{ (1,\dots,1)+2c + 4y\ |\ c\in G_{24}, y\in
\Z^{24}, \text{ and }
\sum_{i=1}^{24} y_i\equiv 1\mod 2\}.
\end{split}
\]

Let $\varepsilon_i= (0,\dots,1,\dots, 0)$ be elements in $\Z^{24}$ such that
the $i$-th entry is $1$ and the other entries are $0$ and set
\[
 \al_{2i-1}=\frac{4}{\sqrt{8}}(\varepsilon_{2i}+\varepsilon_{2i-1}) \quad\text{
and } \quad  \al_{2i}=\frac{4}{\sqrt{8}}(\varepsilon_{2i}-\varepsilon_{2i-1})
\]
for $1\leq i\leq 12$. Then $\{\al_1,\dots,  \al_{24}\}$ forms a 4-frame for
$\Lambda$.

Let $F=\oplus_{i=1}^{24}\Z\al_i$.
Then $\Lambda/F$ determines an extremal type II $\Z_4$-code $\calC$ of
length $24$. The generator matrix of $\calC$ is given  by
\[
\tiny
\left(
\begin{array}{cccccc}
1 1 1 1 & 1 1 1 1 & 1 1 1 1 & 1 1 1 1 & 0 0 0 0 & 0 0 0 0\\
1 1 1 1 & 1 1 1 3 & 2 0 0 0 & 0 0 0 0 & 2 0 0 0 & 0 0 0 0\\
0 0 0 0 & 0 0 0 0 & 1 1 1 1 & 1 1 1 1 & 1 1 1 1 & 1 1 1 1\\
1 1 1 3 & 2 0 0 0 & 1 1 1 1 & 0 0 0 0 & 1 1 1 3 & 0 0 0 0\\
1 3 2 0 & 1 1 0 0 & 1 1 0 0 & 1 1 0 0 & 1 1 0 0 & 1 1 0 0\\
3 2 1 0 & 1 0 1 0 & 1 0 1 0 & 1 0 1 0 & 1 0 1 0 & 1 0 1 0\\
2 2 2 2 & 0 0 0 0 & 0 0 0 0 & 0 0 0 0 & 0 0 0 0 & 0 0 0 0\\
2 2 0 0 & 2 2 0 0 & 0 0 0 0 & 0 0 0 0 & 0 0 0 0 & 0 0 0 0\\
2 0 2 0 & 2 0 2 0 & 0 0 0 0 & 0 0 0 0 & 0 0 0 0 & 0 0 0 0\\
0 0 0 0 & 0 0 0 0 & 2 2 2 2 & 0 0 0 0 & 0 0 0 0 & 0 0 0 0\\
0 0 0 0 & 0 0 0 0 & 2 2 0 0 & 2 2 0 0 & 0 0 0 0 & 0 0 0 0 \\
0 0 0 0 & 0 0 0 0 & 2 0 2 0 & 2 0 2 0 & 0 0 0 0 & 0 0 0 0 \\
2 2 0 0 & 0 0 0 0 & 2 2 0 0 & 0 0 0 0 & 0 0 0 0 & 0 0 0 0\\
2 0 2 0 & 0 0 0 0 & 2 0 2 0 & 0 0 0 0 & 0 0 0 0 & 0 0 0 0\\
2 0 0 0 & 2 0 0 0 & 2 0 0 0 & 2 0 0 0 & 0 0 0 0 & 0 0 0 0\\
0 0 0 0 & 0 0 0 0 & 2 2 0 0 & 0 0 0 0 & 2 2 0 0 & 0 0 0 0\\
0 0 0 0 & 0 0 0 0 & 2 0 2 0 & 0 0 0 0 & 2 0 2 0 & 0 0 0 0\\
0 0 0 0 & 0 0 0 0 & 2 0 0 0 & 2 0 0 0 & 2 0 0 0 & 2 0 0 0
\end {array}
\right)_.
\]

\begin{remark} Harada and Munemasa \cite{HM} checked that the code $\calC$ is the unique extremal type II $\Z_4$-code of length $24$ up to isomorphism such that $\dim \calC_1=6$ and $\calC_0$ has minimum weight $4$.
\end{remark}

The automorphism group of the $\Z_4$-code $\calC$ was computed by Harada and Munemasa using computer \cite{HM}. The automorphism group $\aut(\calC)$ has the shape
\[
2^9.( 2^9. (\Sym_3\times {\rm GL}(3,2))) \quad \text{and} \quad \overline{\aut(\calC)}\cong 2^9.(\Sym_3\times {\rm GL}(3,2)).
\]
Note  also that
$\calC_1=\{ \varphi_2(\al)\ |\ \al\in \calC\} =\Span_{\Z_2}\{
(1^8 0^{16}), (0^81^80^8), (\al,\al,\al)|\, \al\in H_8\},
$
which is a $[24,6,8]$ code (cf. \cite{DGH}). The automorphism group of $\calC_1$ was computed in \cite[Appendix C]{DGH}. The shape is as follows.
\[
 \aut(\calC_1)\cong 2^9.(\Sym_3\times {\rm GL}(3,2)).
\]
Since $\calC$ is self-dual, $\calC_1=\calC_0^\perp$, and $\aut(\calC_1)=\aut(\calC_0)$.

Let $T_{48}$ be the Virasoro frame of the moonshine VOA $V^\natural=\tilde{V}_\Lambda$ associated
to $F$. Then the structure codes $(C,D)$ of $V^\natural$ associated
to $T_{48}$ are given by $C=D^\perp $ and
\[\small
 D=
\left(
\begin{matrix}
1111\ 1111 \ 1111 \ 1111 \ 1111\ 1111 \ 1111 \ 1111 \ 1111 \ 1111 \ 1111 \
1111\\
1111\ 1111 \ 1111 \ 1111 \ 1111\ 1111 \ 1111 \ 1111 \ 0000 \ 0000 \ 0000 \
0000\\
1111\ 1111 \ 1111 \ 1111 \ 0000 \ 0000 \ 0000 \
0000 \ 0000 \ 0000 \ 0000 \ 0000\\
1111\ 1111  \ 0000 \ 0000 \ 1111 \ 1111 \ 0000 \
0000 \ 1111 \ 1111 \ 0000 \ 0000\\
1111\ 0000 \ 1111 \ 0000 \ 1111 \ 0000 \ 1111  \ 0000 \ 1111 \ 0000 \ 1111  \ 0000\\
1100\ 1100 \ 1100 \ 1100 \ 1100 \ 1100\ 1100\ 1100 \ 1100 \ 1100 \ 1100 \ 1100\\
1010\ 1010 \ 1010 \ 1010 \ 1010 \ 1010\ 1010\ 1010 \ 1010 \ 1010 \ 1010 \ 1010\\
\end{matrix}
\right)_.
\]
Note that $D=\Span_{\Z_2}\{(1^{16}, 0^{16}, 0^{16}), (0^{16}, 1^{16},  0^{16}), (\alpha,\alpha,\alpha) \mid
  \ \alpha\in \RM(1,4)\}$.  In this case,
\[
\aut(C)=\aut(D)\cong 2^{12}.(\Sym_3\times{\rm {\rm GL}}(4,2))
\]
(see \cite[Appendix C]{DGH} for details). Moreover,
\[
\begin{split}
P&= \{ \gamma \in \Z_2^{48} \mid \alpha\cd \gamma \in C\ \text{for all}\
  \alpha\in {D}\}\\
  &  = \{ (\alpha,\beta,\gamma) \in \Z_2^{48} \mid \alpha,\beta,\gamma\in \RM(2,4)
  \ \text{and}\ \alpha+\beta+\gamma\in \RM(1,4)\}
\end{split}
\]
and $\dim P= 11+11+5=27$. Thus, $\pstab_{\aut(V^\natural)}(T_{48})$ has the shape $2^{7+20}$.

Since $|\aut(C): {K}|=|\aut(\calC_0): \overline{\aut(\calC)}|=1$, we have $$\stab_{\Aut(V^\natural)}(T_{48})/\pstab_{\Aut(V^\natural)}(T_{48})\cong \aut(C),$$ and the frame stabilizer
$\stab_{\Aut(V^\natural)}(T_{48})$ has the shape $2^{7+20}. ( 2^{12}.(\Sym_3\times{\rm GL}(4,2)))$.

\begin{remark}
It was also shown in \cite[Lemma 9.3]{M3} that $$\stab_{\Aut(V^\natural)}(T_{48})/\pstab_{\Aut(V^\natural)}(T_{48})\cong \aut(C).$$
\end{remark}

\end{document}